\documentclass[11pt]{amsart}

\input epsf

\usepackage{amsmath}
\usepackage{amssymb}
\usepackage{amsfonts}
\usepackage[mathscr]{eucal}
\usepackage{graphicx}
\usepackage{epstopdf}
\usepackage{verbatim}
\usepackage{amscd}
\usepackage{upgreek}
\usepackage{dsfont}
\usepackage{caption}
\usepackage{amsthm}
\usepackage{color}
\usepackage{hyperref}
\usepackage[all]{xy}
\usepackage{pdfsync}

\usepackage{cleveref}
\crefdefaultlabelformat{#2{\scshape #1}#3}
\crefname{subfigure}{{\scshape Fig.}}{{\scshape Figs.}}
\Crefname{subfigure}{{\scshape Figure}}{{\scshape Figures}}
\crefname{figure}{{\scshape fig.}}{{\scshape figs.}}
\Crefname{figure}{{\scshape Figure}}{{\scshape Figures}}

\usepackage{subcaption}

\usepackage{ccfonts,eulervm}
\usepackage[T1]{fontenc}

\makeindex
\numberwithin{equation}{section}

\theoremstyle{definition}
\newtheorem*{Definition}{Definition}

\usepackage{geometry}  
\usepackage[parfill]{parskip}  

\usepackage[percent]{overpic}   

\title{Ma-Schlenker \textit{c}-Octahedra in the $2$-Sphere}

\author{John C. Bowers}
\address{Department of Computer Science, James Madison University, 
Harrisonburg VA 22807}
\email{johnchristopherbowers@gmail.com}
\author[Philip L. Bowers]{Philip L. Bowers}
\address{Department of Mathematics, The Florida State University, 
Tallahassee FL 32306}
\email{bowers@math.fsu.edu}

\date{\today} 

\begin{document}

\begin{abstract}
We present constructions inspired by the Ma-Schlenker example of~\cite{Ma:2012hl} that show the non-rigidity of spherical inversive distance circle packings. In contrast to the use in~\cite{Ma:2012hl} of an infinitesimally flexible Euclidean polyhedron, embeddings in de Sitter space, and Pogorelov maps, our elementary constructions use only the inversive geometry of the $2$-sphere.
\end{abstract}

\maketitle

\section*{Introduction}
In~\cite{Bowers:2004bg}, P.~Bowers and K.~Stephenson questioned whether inversive distance circle packings of surfaces are uniquely determined by their underlying triangulation and the inversive distances between the pairs of adjacent circles. Guo~\cite{Guo:2011kf} confirmed the local rigidity of inversive distance circle packings on closed orientable surfaces of non-negative genus, and subsequently Luo~\cite{Luo:2011ex} verified the global rigidity of these packings, answering the Bowers-Stephenson question in the affirmative. Contrasted to this is the beautiful and surprising example of Ma and Schlenker in~\cite{Ma:2012hl} that provides a counterexample to uniqueness in the spherical case. They produced pairs of packing radii that determine pairs of geodesic triangulations and circle packings on the $2$-sphere $\mathbb{S}^{2}$ realizing the same inversive distance data, but for which there is no inversive transformation taking one of the circle patterns to the other. This was doubly surprising as the famous Koebe-Andre'ev-Thurston Circle Packing Theorem implies uniqueness of spherical packings up to inversive equivalence whenever the edge labels are all in the unit interval---the case of tangent or overlapping circle packings.

The ingredients of Ma and Schlenker's example are Sch\"onhardt's twisted octahedron (an infinitesimally flexible polyhedron in Euclidean space $\mathbb{E}^{3}$), embeddings in de Sitter space $\mathbb{S}^{3}_{1}$, and special properties of the Pogorelov map between different geometries. In contrast, we provide a construction of a large family of Ma-Schlenker-like examples using only inversive geometry. In fact, we show how to construct many counterexamples to the uniqueness of inversive distance circle packings in the $2$-sphere.

\begin{figure}
\begin{subfigure}[b]{0.3\textwidth}
\includegraphics[width=\textwidth]{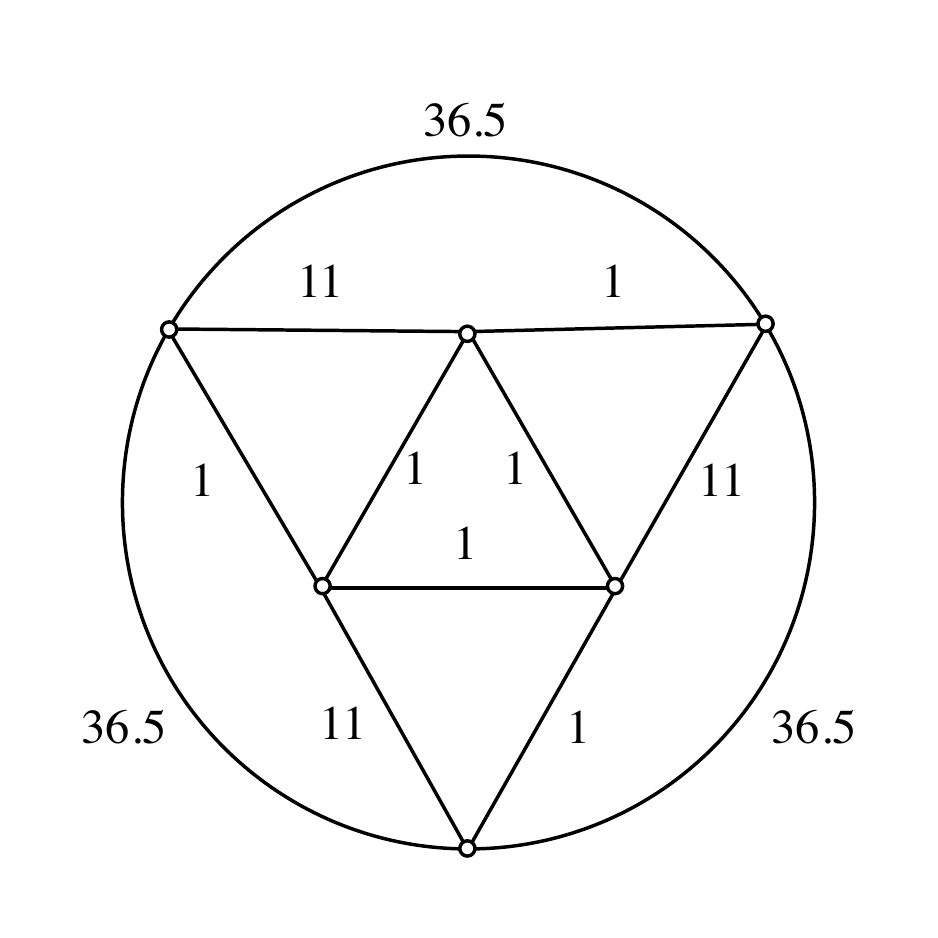}
\caption{Octahedral graph labeled with the inversive distances from \cref{fig:criticalb}.}
\label{fig:criticala}
\end{subfigure}
\quad
\begin{subfigure}[b]{0.3\textwidth}
\includegraphics[width=\textwidth]{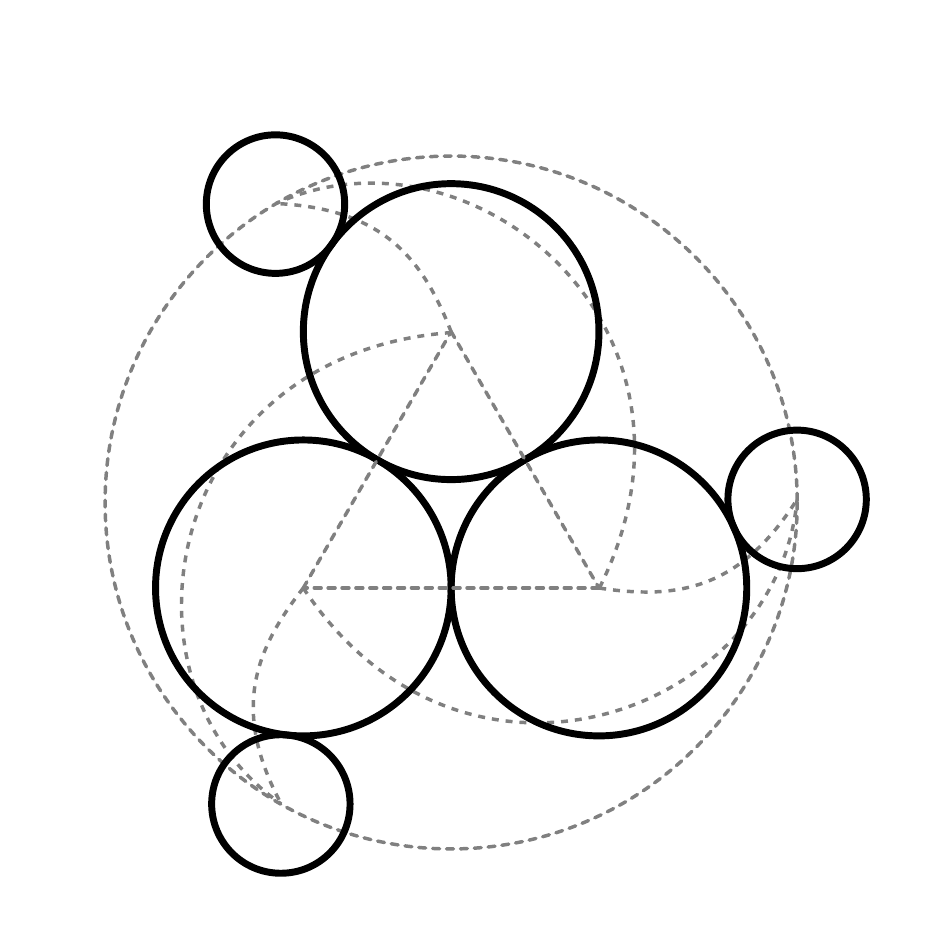}
\caption{A planar circle pattern realizing the octahedral graph in \cref{fig:criticala}.}
\label{fig:criticalb}
\end{subfigure}
\quad
\begin{subfigure}[b]{0.3\textwidth}
\includegraphics[width=\textwidth]{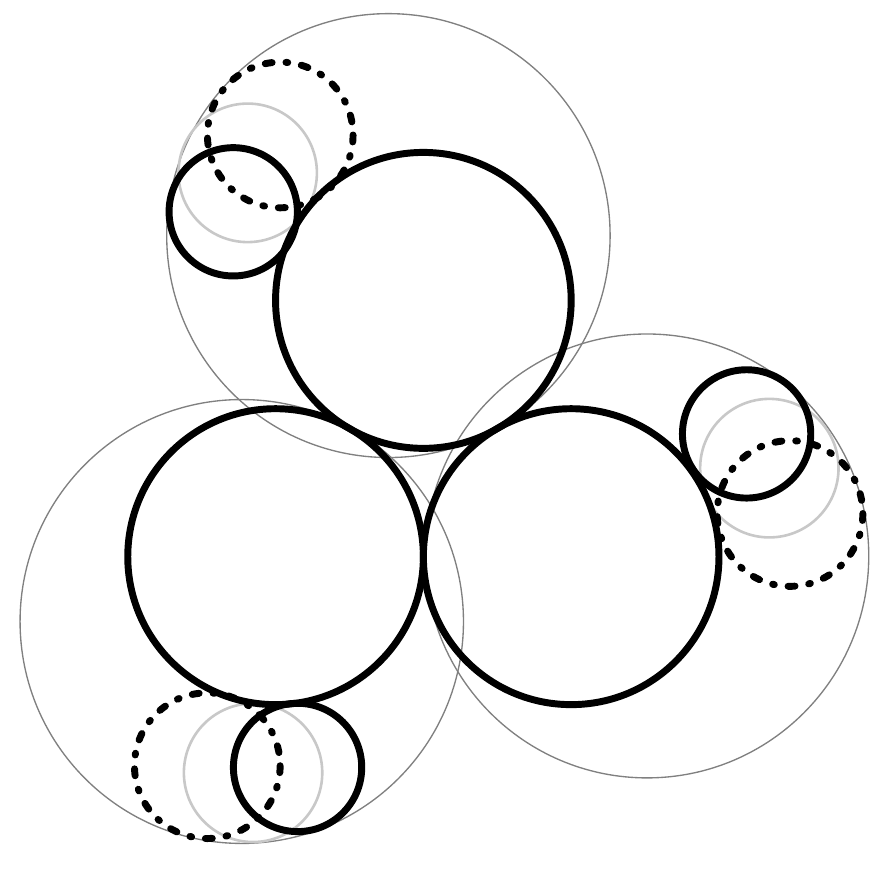}
\caption{Two non-equivalent realizations with inversive distance 37 on the outer edges.}
\label{fig:criticalc}
\end{subfigure}
\caption{A critical Ma-Schlenker circle octahedron and two non-inversive equivalent nearby patterns with the same inversive distances on the graph.}\label{fig:critical}
\end{figure}
\Cref{fig:criticalb} pictures a planar configuration of six circles, three inner of a common radius and three outer of a common radius. Edges between selected circle centers representing inversive distance constraints are depicted by dashed gray lines. The edges have the incidence relations of the $1$-skeleton of an octahedron and by stereographic projection to the $2$-sphere $\mathbb{S}^{2}$ followed by an appropriate M\"obius transformation, may be arranged so that the labeled inversive distances are preserved, and so that connecting centers of adjacent circles by great circle arcs cuts out an octahedral triangulation of $\mathbb{S}^{2}$. \Cref{fig:criticala} depicts the underlying octahedral graph labeled with the inversive distances between the adjacent circles from \cref{fig:criticalb}, accurate to the nearest tenth.  This determines a \textit{critical Ma-Schlenker circle octahedron}, an inversive distance circle packing of the $2$-sphere with the following property. By varying the three outer circles as pictured in \cref{fig:criticalc}, one may obtain infinitely many pairs of circle packings such that the two circle packings of a pair are not inversive equivalent even though they have the same inversive distances labeling corresponding edges. This was the first of many similar examples the authors constructed, and the remainder of this paper develops the tools needed to construct such examples and verify the claims of this paragraph.

\section{Preliminaries from Circle Packing and Inversive Geometry}
In this preliminary section, we recall and expand upon the basic facts about inversive distance, circle packings, and inversive geometry that are needed to describe and understand our examples. 
\subsection{Inversive distance in the plane and the $2$-sphere}
Let $C_{1}$ and $C_{2}$ be distinct circles in the complex plane $\mathbb{C}$ centered at the respective points $p_{1}$ and $p_{2}$, of respective radii $r_{1}$ and $r_{2}$, and bounding the respective \textit{companion disks} $D_{1}$ and $D_{2}$. The \textit{inversive distance} $\langle C_{1}, C_{2} \rangle$ between $C_{1}$ and $C_{2}$ is
\begin{equation}\label{EQ:planeID} 
	\langle C_{1}, C_{2} \rangle = \frac{|p_{1}-p_{2}|^{2} - r_{1}^{2} - r_{2}^{2}}{2r_{1}r_{2}}.
\end{equation}
We find it convenient to extend this to disk pairs by defining the \textit{inversive distance} $\langle D_{1},D_{2}\rangle$ between the disks $D_{1}$ and $D_{2}$ in exactly the same way, with $\langle D_{1}, D_{2} \rangle = \langle C_{1},C_{2}\rangle$. The \textit{absolute inversive distance} between distinct circles is the absolute value of the inversive distance and is a M\"obius invariant of the placement of two circles in the plane. This means that there is a M\"obius transformation of $\mathbb{C}$ taking one circle pair to another if and only if the absolute inversive distances of the two pairs agree.\footnote{There is a similar statement for the inversive distance without the modifier \textit{absolute}. Indeed, the inversive distance is a M\"obius invariant of the placement of two relatively oriented circles in the plane. See Bowers-Hurdal~\cite{Bowers:2003kr}.} The important geometric facts that make the inversive distance useful in inversive geometry and circle packing are as follows. When $\langle C_{1},C_{2}\rangle > 1$, $D_{1} \cap D_{2} = \emptyset$ and $\langle C_{1},C_{2}\rangle = \cosh \delta$, where $\delta$ is the hyperbolic distance between the totally geodesic hyperbolic planes in the upper-half-space model $\mathbb{C} \times (0,\infty)$ of $\mathbb{H}^{3}$ whose ideal boundaries are $C_{1}$ and $C_{2}$. When $\langle C_{1},C_{2}\rangle = 1$, $D_{1}$ and $D_{2}$ are tangent at their single point of intersection. When $1 >\langle C_{1},C_{2}\rangle \geq 0$, $D_{1}$ and $D_{2}$ overlap with angle $0 < \theta \leq \pi/2$ with $\langle C_{1},C_{2}\rangle = \cos \theta$. In particular, $\langle C_{1},C_{2}\rangle = 0$ precisely when $\theta = \pi/2$. When $\langle C_{1},C_{2}\rangle < 0$, then $D_{1}$ and $D_{2}$ overlap by an angle greater than $\pi/2$. This includes the case where one of $D_{1}$ or $D_{2}$ is contained in the other, this when $\langle C_{1},C_{2}\rangle \leq -1$. In fact, when $\langle C_{1},C_{2}\rangle < -1$ then $\langle C_{1},C_{2}\rangle = - \cosh \delta$ where $\delta$ has the same meaning as above, and when $\langle C_{1},C_{2}\rangle = -1$ then $C_{1}$ and $C_{2}$ are `internally' tangent. When $-1 < \langle C_{1},C_{2}\rangle <0$, then the overlap angle of $D_{1}$ and $D_{2}$ satisfies $\pi > \theta > \pi/2$ and again $\langle C_{1},C_{2}\rangle = \cos \theta$.

In the $2$-sphere $\mathbb{S}^{2}$, the inversive distance may be expressed as\footnote{In both Luo~\cite{Luo:2011ex} and Ma-Schlenker~\cite{Ma:2012hl} there is a typo in the expression for the spherical formula for inversive distance. They report the negative of this formula. A quick 2nd order Taylor approximation shows that this formula reduces to Expression~\ref{EQ:planeID} in the limit as the arguments of the sines and cosines approach zero.} 
\begin{equation}\label{EQ:sphereID} 
	\langle C_{1}, C_{2} \rangle = \frac{-\cos \sphericalangle ( p_{1},p_{2} ) + \cos(r_{1}) \cos(r_{2})}{\sin(r_{1}) \sin(r_{2})}.
\end{equation}
Here, $\sphericalangle (p_{1},p_{2})$ denotes the spherical distance between the centers, $p_{1}$ and $p_{2}$, of the respective circles\footnote{Any circle in $\mathbb{S}^{2}$ bounds two distinct disks. Without explicitly stating so, we always assume that one of these has been chosen as a companion disk. The center and radius of a circle in $\mathbb{S}^{2}$ are the center and radius of its companion disk. The ambiguity should cause no confusion.} $C_{1}$ and $C_{2}$ with respective spherical radii $r_{1}$ and $r_{2}$.\footnote{It is in no way obvious that Formul{\ae}~\ref{EQ:planeID} and \ref{EQ:sphereID} are M\"obius invariants of circle pairs. There is a not so well-known development of inversive distance, which applies equally in spherical, Euclidean, and hyperbolic geometry, that uses the cross-ratio of the four points of intersection of $C_{1}$ and $C_{2}$ with a common orthogonal circle. It is computationally less friendly than~\ref{EQ:planeID} and \ref{EQ:sphereID}, but has the theoretical advantage of being manifestly M\"obius-invariant. See~\cite{Bowers:2003kr}.} Stereographic projection to the plane $\mathbb{C}$ preserves the absolute inversive distance of circle pairs and, as long as neither of the two companion disks $D_{1}$ and $D_{2}$ contains the north pole, it preserves the inversive distance.

\subsection{Edge-labeled triangulations of the $2$-sphere and inversive distance circle packings}
Bowers and Stephenson originally introduced inversive distance circle packings in~\cite{Bowers:2003kr} with inversive distances restricted to be non-negative. We will take this opportunity to define packings with no restrictions on the inversive distance and offer some warnings of the pitfalls of the more general setting. We are concerned here with configurations of circles in the $2$-sphere $\mathbb{S}^{2}$ with a specified pattern of inversive distances. Let $K$ be an abstract oriented triangulation of $\mathbb{S}^{2}$ and $\beta : E(K) \to \mathbb{R}$ a mapping defined on $E(K)$, the set of edges of $K$. Call $K$ together with $\beta$ an \textit{edge-labeled triangulation} with \textit{edge label} $\beta$ and denote it as $K_{\beta}$. We denote an edge of $K$ with vertices $u$ and $v$ by $uv$, and an oriented face with vertices $u$, $v$, and $w$ ordered respecting the orientation of $K$ by $uvw$. We define two types of circle configurations that realize the inversive distance data encoded in an edge-labeled triangulation.
\begin{Definition}
	A \textit{circle realization} for $K_{\beta}$ is a collection $\mathscr{C} = \{C_{v} : v \in V(K)\}$ of circles in either the plane $\mathbb{C}$ or the $2 $-sphere $\mathbb{S}^{2}$ indexed by the vertex set $V(K)$ of $K$ such that $\langle C_{u}, C_{v}\rangle = \beta (uv )$ whenever $uv$ is an edge of $K$. When $uv$ is an edge of $K$, the corresponding circles $C_{u}$ and $C_{v}$ are said to be \textit{adjacent}.
\end{Definition}
A circle packing for $K_{\beta}$ is a circle realization in the $2$-sphere where the circles are placed in a way that the circles $C_{u}$, $C_{v}$, and $C_{w}$ form a positively oriented triple in $\mathbb{S}^{2}$ whenever $uvw$ is a positively oriented face of $K$. This general definition allows for behavior, for example branch structures, that we wish to avoid. Moreover, this general definition involves some subtleties that arise from the analytic fact that there is no M\"obius-invariant metric on the $2$-sphere, and from the topological fact that a simple closed curve in the $2$-sphere fails to have a well-defined inside as it has two complementary domains, both of which are topological disks. For one example of these subtleties, centers and radii of circles are not well-defined in inversive geometry so that two circle realizations for $K_{\beta}$ may have differing properties, say with one forming a triangulation of $\mathbb{S}^{2}$ by connecting adjacent centers along great circular arcs and the other not, even though one is the M\"obius image of the other; see~\cite{Bowers:rigidityOfCircle:techReport} for enlightening discussion and examples.\footnote{This behavior does not occur for the traditional tangency and overlapping packings.} We choose a more restrictive definition that conforms to the construction of circle packings by the use of polyhedral metrics, a construction we review below and that is used by Ma and Schlenker in their construction. In particular, our interest is in circle packings that produce  isomorphic copies of the triangulation $K$ by connecting the centers of adjacent circles by geodesic segments. In this paper then, circle packing means the following.
\begin{Definition} An \textit{inversive distance circle packing}, or simply a \textit{circle packing}, for $K_{\beta}$ is a collection $\mathscr{C} = \{C_{v} : v \in V(K)\}$ of circles in $\mathbb{S}^{2}$ with four properties:
\begin{enumerate}
\item[(i)] $\mathscr{C}$ is a circle realization for $K_{\beta}$;
\item[(ii)] when $uv$ is an edge of $K$, the centers of $C_{u}$ and $C_{v}$ are not antipodal;\footnote{Though (iii) implies (ii), we opt to state (ii) explicitly, for emphasis.}
\item[(iii)] when $uvw$ is a face of $K$, the centers of $C_{u}$, $C_{v}$, and $C_{w}$ do not lie on a great circle of $\mathbb{S}^{2}$;
\item[(iv)] joining all the pairs of centers of adjacent circles $C_{u}$ and $C_{v}$ by geodesic segments of $\mathbb{S}^{2}$ produces a triangulation of $\mathbb{S}^{2}$, necessarily isomorphic with $K$.
\end{enumerate}
\end{Definition}
These conditions imply that when $uvw$ is a face of $K$, the centers of $C_{u}$, $C_{v}$, and $C_{w}$ do not lie on a common geodesic, and for two distinct faces of $K$, the interiors of the corresponding geodesic triangles determined by $\mathscr{C}$ have empty intersection.\footnote{All of this generalizes in a straightforward way to triangulations of arbitrary constant curvature surfaces, closed or not, but our concern will be with the $2$-sphere. Also, we have described here the univalent circle packings, so this discussion can be generalized to packings that are only locally univalent and even to ones with branch vertices. See~\cite{kS05} and \cite{Bowers:1996ur}.} To reiterate the warning stated in the preceding paragraph, under this restricted definition of circle packing, the M\"obius image of a circle packing need not be a circle packing.

 When $\mathscr{C}$ is a circle packing for $K_{\beta}$, the set of radii $r_{v}$ of the circles $C_{v} \in \mathscr{C}$ is a set of \textit{packing radii} for $K_{\beta}$. Two circle packings $\mathscr{C}$ and $\mathscr{C}'$ for $K_{\beta}$ are \textit{inversive equivalent} if there is an inversive transformation $T\in \mathrm{Inv}(\mathbb{S}^{2})$ for which $T(C_{v}) =C_{v}'$ for each vertex $v$. Here of course $\mathrm{Inv}(\mathbb{S}^{2})$ is the group of inversive transformations of the $2$-sphere generated by inversions through the circles of $\mathbb{S}^{2}$. The packings are \textit{M\"obius equivalent} if $T$ can be chosen to be a M\"obius transformation, an element of $\text{M\"ob}(\mathbb{S}^{2})$, the M\"obius group of the $2$-sphere generated by even numbers of compositions of inversions. The corresponding sets of packing radii, $\{r_{v} : v \in V(K)\}$ and $\{r_{v}': v\in V(K)\}$, are said to be \textit{inversive equivalent} provided $\mathscr{C}$ and $\mathscr{C}'$ are inversive equivalent. The question of uniqueness of packings or of packing radii for $K_{\beta}$ is always up to inversive equivalence of packings or packing radii, and our interest is in constructing pairs of circle packings for the same edge-labeled triangulation that are not inversive equivalent.

\begin{figure}
\includegraphics[width=0.7\textwidth]{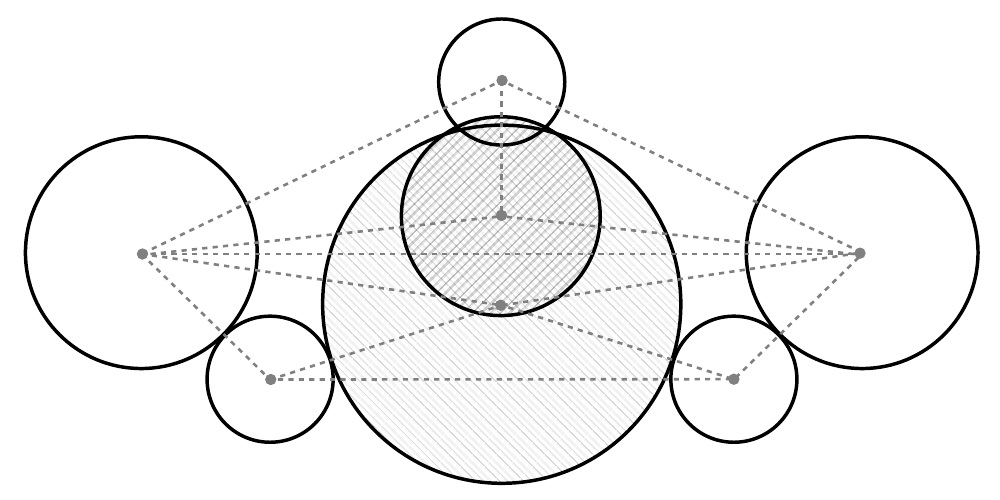}
\caption{An edge-segregated inversive-distance packing that is not globally segregated. The overlap angle for the two shaded circles is greater than $\pi / 2$, but there is no edge between them.}
\label{fig:locallysegregated}
\end{figure}
Circle packings traditionally have been studied when adjacent circles overlap non-trivially, and with angles of overlap at most $\pi/2$; i.e., when $\beta$ takes values in the unit interval $[0,1]$.\footnote{But see Rivin~\cite{Rivin:1994kg}.} This is because of various geometric, computational, and theoretical difficulties associated to overlaps greater than $\pi/2$. The Koebe-Andre'ev-Thurston Circle Packing Theorem applies only to these, as do many of the known existence and uniqueness results.\footnote{In fact, the global rigidity result of Luo~\cite{Luo:2011ex} proving uniqueness for packing radii on closed parabolic and hyperbolic surfaces holds only with the assumption of overlaps of at most $\pi/2$. The authors in~\cite{Bowers:rigidityOfCircle:techReport} study the rigidity of \textit{circle frameworks}, which are more general collections of circles than are circle realizations and packings, with no bounds on the overlaps of the companion disks. There we show that the Luo rigidity results fail to hold when circle packings with overlaps greater than $\pi/2$ are allowed.} The Bowers-Stephenson question of the uniqueness of inversive distance circle packings was asked for packings, still with overlaps at most $\pi/2$, but where adjacent circles may be \textit{separated}, when $\beta$ may take values greater than unity. This motivates the following definitions. A collection of circles in the plane or the sphere is said to be \textit{segregated} provided the corresponding collection of companion disks bounded by the circles pairwise overlap by no more than $\pi/2$. A circle realization or packing $\mathscr{C}$ for $K_{\beta}$ is said to be \textit{edge-segregated} if the companion disks of any two adjacent circles overlap by at most $\pi/2$, i.e., if $\beta$ takes on only non-negative values.  A segregated circle packing obviously is edge-segregated, but an edge-segregated circle packing need not be segregated; see \cref{fig:locallysegregated}. A circle packing that is not edge-segregated is said to have \textit{deep overlaps} at the adjacent circle pairs where $\beta$ is negative. Finally, a collection of circles is \textit{separated} provided the corresponding collection of companion disks are pairwise disjoint, and a circle realization or packing is \textit{edge-separated} if $\beta$ only takes values greater than unity so that adjacent circles are separated.

The Koebe-Andre'ev-Thurston Circle Packing Theorem is the fundamental existence and uniqueness result for circle packings for which the image of $\beta$ is contained in the unit interval. There are no good results on the existence of general circle packings, edge-segregated or not, when $\beta$ takes values outside the unit interval $[0,1]$.  The Guo, Luo, and Ma-Schlenker results already referenced are the current state of the art in uniqueness results for edge-segregated circle packings. The reference~\cite{Bowers:rigidityOfCircle:techReport} makes some observations on rigidity for general circle realizations and packings, where $\beta$ is completely unrestricted, and shows the non-uniqueness of packing radii for general packings of arbitrary surfaces as well as the non-rigidity of segregated realizations.

An alternate way to describe circle packings in the $2$-sphere with prescribed inversive distances is through the constructive use of spherical polyhedral surfaces. Let $\ell: E(K) \to (0,\infty)$ be a map that satisfies the strict triangle inequality for the edges of any face $uvw$, meaning that 
\begin{equation}\label{EQ:StrictTI}
	\ell(uw)  < \ell(uv)  + \ell(vw).
\end{equation}
If in addition
\begin{equation}\label{EQ:GreatCircle}
	\ell(uv) + \ell(vw)  + \ell(wu) < 2\pi
\end{equation} 
for each face $uvw$, we call $\ell$ a \textit{spherical length function} for $K$. Associated to a triangulation $K$ and a spherical length function $\ell$ is the \textit{spherical polyhedral surface $S = S(K, \ell)$} obtained by gluing together in the pattern of $K$ spherical triangles whose side lengths are given by $\ell$. Each edge $uv$ of $K$ is identified metrically as a Euclidean segment of length $\ell (uv)$ and each face $uvw$ of $K$ is identified metrically with a spherical triangle in $\mathbb{S}^{2}$ whose side-lengths are $\ell(uv)$, $\ell(vw)$, and $\ell(wu)$. The surface $S(K,\ell)$ topologically is a $2$-sphere with a singular Riemannian metric of constant curvature $+1$. Any singularities occur at the vertices of $K$, and this only when the angle sum of the spherical triangles that meet at the vertex is other than $2\pi$. 

Associated to any circle packing $\mathscr{C}$ of an edge-labeled triangulation $K_{\beta}$ is an obvious spherical length function $\ell$ defined by letting $\ell (uv)$ be the spherical distance in $\mathbb{S}^{2}$ between the centers of adjacent circles $C_{u}$ and $C_{v}$. The corresponding spherical polyhedral surface $S(K,\ell)$ then is isometric with $\mathbb{S}^{2}$. Now this spherical length function $\ell$ can be described completely in terms of the edge label $\beta$ and the set of packing radii $\{r_{v}: v\in V(K)\}$ using Formula~\ref{EQ:sphereID}. This hints at a possible approach to finding inversive distance circle packings for a given edge-labeled triangulation $K_{\beta}$. Starting with $K_{\beta}$ and a set of positive numbers $0 < r_{v} < \pi$ for $v\in V(K)$, define the function $\ell: E(K) \to (0, \infty)$ by
\begin{equation}
	\ell (uv) = \cos^{-1} \left(\cos r_{u} \cos r_{v} - \beta(uv) \sin r_{u} \sin r_{v} \right),
\end{equation}
provided that $|\cos r_{u} \cos r_{v} - \beta(uv) \sin r_{u} \sin r_{v}| < 1$ for all edges $uv$ of $K$. If the \textit{proposed packing radii} $r_{v}$ are chosen so that $\ell$ exists and satisfies Inequalities~\ref{EQ:StrictTI} and \ref{EQ:GreatCircle} for each face $uvw$, then $\ell$ is a spherical length function for $K$ and the spherical polyhedral surface $S(K, \ell)$ supports a collection of metric circles $\mathscr{C} =\{C_{v} : v \in V(K)\}$, where $C_{v}$ is centered at the vertex $v$  with radius $r_{v}$. Though $S(k, \ell)$ is singular for typical proposed packing radii, the idea is to vary the proposed radii $r_{v}$ in the hopes of finding a set that removes all the singularities. When this occurs, the surface $S(K,  \ell)$ is isometric to the standard $2$-sphere and the collection $\mathscr{C}$ is a circle packing for $K_{\beta}$, with $\langle C_{u},C_{v} \rangle = \beta (uv)$ for every edge $uv$ of $K$. Of course there are a lot of `ifs' here, and many edge-labeled triangulations will not have any circle packings. The general existence question is quite intricate and is not pursued here as our interest is rigidity.

\subsection{M\"obius flows}
The most important ingredients from inversive geometry needed to construct our examples are the M\"obius flows associated to two distinct circles in the extended complex plane. Here are the pertinent facts. Circles $C_{1} \neq C_{2}$ lie in a unique \textit{coaxial family} $\mathcal{A}_{C_{1},C_{2}}$ of circles in the extended complex plane $\widehat{\mathbb{C}}$ whose elements serve as the flow lines of certain $1$-parameter subgroups of the M\"obius group $\text{M\"ob}(\widehat{\mathbb{C}})\cong \mathrm{PSL}(2, \mathbb{C})$ acting on the extended plane as linear fractional transformations. The family $\mathcal{A}_{C_{1},C_{2}}$ is invariant under these flows. When $C_{1}$ and $C_{2}$ are disjoint, any such flow is an \textit{elliptic flow} conjugate in $\text{M\"ob}(\widehat{\mathbb{C}})$ to a standard rotation flow of the form $t \mapsto R_{\lambda t}$, where $\lambda \neq 0$ and $R_{\lambda t}$ is the rotation $z \mapsto e^{\lambda \boldsymbol{i} t} z$. When $C_{1}$ and $C_{2}$ meet in a single point $p$, the flow is \textit{parabolic} and is conjugate in $\text{M\"ob}(\widehat{\mathbb{C}})$ to a standard translation flow of the form $t\mapsto T_{\lambda t}$, where $\lambda \neq 0$ and $T_{\lambda t}$ is the translation $z \mapsto z + \lambda t$. Finally, when $C_{1}$ and $C_{2}$ meet in two distinct points $a$ and $b$, the flow is \textit{hyperbolic} and is conjugate in $\text{M\"ob}(\widehat{\mathbb{C}})$ to a standard scaling flow of the form $t \mapsto S_{\lambda t}$, where $\lambda \neq 0$ and $S_{\lambda t}$ is the scaling map $z \mapsto e^{\lambda t} z$. All the flows determined by two fixed circles $C_{1}$ and $C_{2}$ are said to be \textit{equivalent flows}, with two such flows differing only in the value of the parameter $\lambda$, with $|\lambda|$ the \textit{speed} of the flow. Note that any two distinct circles in the coaxial family $\mathcal{A}_{C_{1},C_{2}}$ determine the same M\"obius flows as $C_{1}$ and $C_{2}$. Of course, any M\"obius flow preserves the inversive distances between circles.

Each coaxial family $\mathcal{A}$ has an associated \textit{orthogonal complement} $\mathcal{A}^{\perp}$, this also a coaxial family for which each circle of $\mathcal{A}$ is orthogonal to each circle of $\mathcal{A}^{\perp}$. In fact, $\mathcal{A}^{\perp}$ is exactly the collection of circles (and lines) that are orthogonal to every member of $\mathcal{A}$, and of course $\mathcal{A}^{\perp\perp} = \mathcal{A}$. Any flow whose flow lines are the circles of $\mathcal{A}$ is generated by inversions through the circles of the orthogonal complement $\mathcal{A}^{\perp}$. Other than the three standard cases of a family of concentric circles (flow lines of a standard rotation flow), of parallel lines (flow lines of a standard translation flow), and of a pencil of lines through a fixed point (flow lines of a standard scaling flow), each coaxial family $\mathcal{A}$ in the plane has a unique line (rather than circle) among its members, this called the \textit{radical axis} of $\mathcal{A}$. The centers of all the circles of $\mathcal{A}$ lie on a common line, the \textit{line of centers} of $\mathcal{A}$. Beautifully, the radical axis of $\mathcal{A}$ is the line of centers of $\mathcal{A}^{\perp}$, and vice-versa. An important property of the pair $\{\mathcal{A}, \mathcal{A}^{\perp}\}$ in the plane is that every circle $C$ centered on the radical axis of $\mathcal{A}$ and orthogonal to a single member of $\mathcal{A}$ other than its radical axis is orthogonal to every member of $\mathcal{A}$, which in turn implies that $C\in \mathcal{A}^{\perp}$. This is just a special case of the fact that any circle orthogonal to both $C_{1}$ and $C_{2}$ is necessarily a member of the orthogonal complement $\mathcal{A}_{C_{1},C_{2}}^{\perp}$.\footnote{For a nice treatment of coaxial families in the plane, see~\cite{Brannan:2011td}, and for a broader treatment that develops both an intrinsic and extrinsic version for the $2$-sphere, see~\cite{Bowers:isoInversive:techReport}.}

\section{Ma-Schlenker \textit{c}-Octahedra---the Examples writ Large}\label{Section:MAOWL}
In this section we detail examples of pairs of spherical circle packings for a fixed edge-labeled triangulation that fail to be inversive equivalent.  These were discovered by studying the properties of the first family of such examples constructed by Ma and Schlenker in~\cite{Ma:2012hl} using fairly sophisticated geometric constructions. Ours are constructed using M\"obius flows and, having properties reminiscent of the Ma-Schlenker examples, will be named after them. This section presents a general construction of such examples and describes their important properties, and the next verifies the claimed properties.

\begin{figure}
\includegraphics[width=0.35\textwidth]{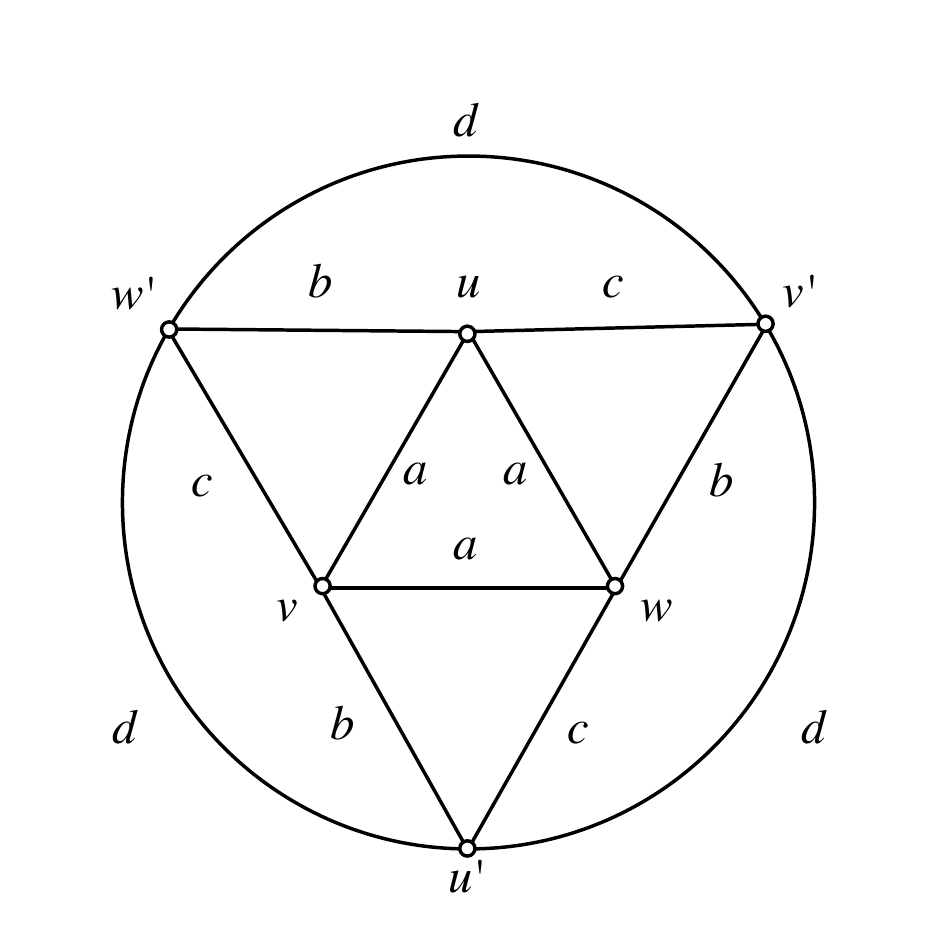}
\caption{The typical Ma-Schlenker octahedron $\mathscr{O}(a,b,c,d)$.}
\label{fig:msocta}
\end{figure}

Let $\mathscr{O}$ be the octahedral triangulation of the $2$-sphere with six vertices, each of valence four, twelve edges, and eight faces, and combinatorially equivalent to the boundary of a regular octahedron. The $1$-skeleton graph of $\mathscr{O}$ is shown embedded in the plane in \cref{fig:criticala}. 
\begin{Definition}
	The edge-labeled triangulation $\mathscr{O}_{\beta}$ is called a \textit{Ma-Schlenker octahedron} provided $\beta$ takes a constant value $a\geq 0$ on the edges of a fixed face $uvw$, a constant value $d\geq 0$ on the edges of its opposite face $w'v'u'$, and alternates between the values $b\geq 0$ and $c\geq 0$ on the `teepee' of edges connecting the vertices of these two opposite faces, as in \cref{fig:msocta}. When we need to emphasize the values of $a$, $b$, $c$, and $d$, we denote $\mathscr{O}_{\beta}$ as $\mathscr{O}(a,b,c,d)$.
\end{Definition}
\begin{Definition}
A circle packing for a Ma-Schlenker octahedron is called a \textit{Ma-Schlenker circle-octahedron}, or a \textit{Ma-Schlenker \textit{c}-octahedron} for short. A pair of circle packings, $\mathscr{C}$ and $\mathscr{C}'$, for the same Ma-Schlenker octahedron $\mathscr{O}_{\beta}$ is called a \textit{Ma-Schlenker pair} provided $\mathscr{C}$ and $\mathscr{C}'$ are not inversive equivalent.
\end{Definition}
\begin{figure}
\begin{subfigure}[b]{0.48\textwidth}
\includegraphics[width=\textwidth]{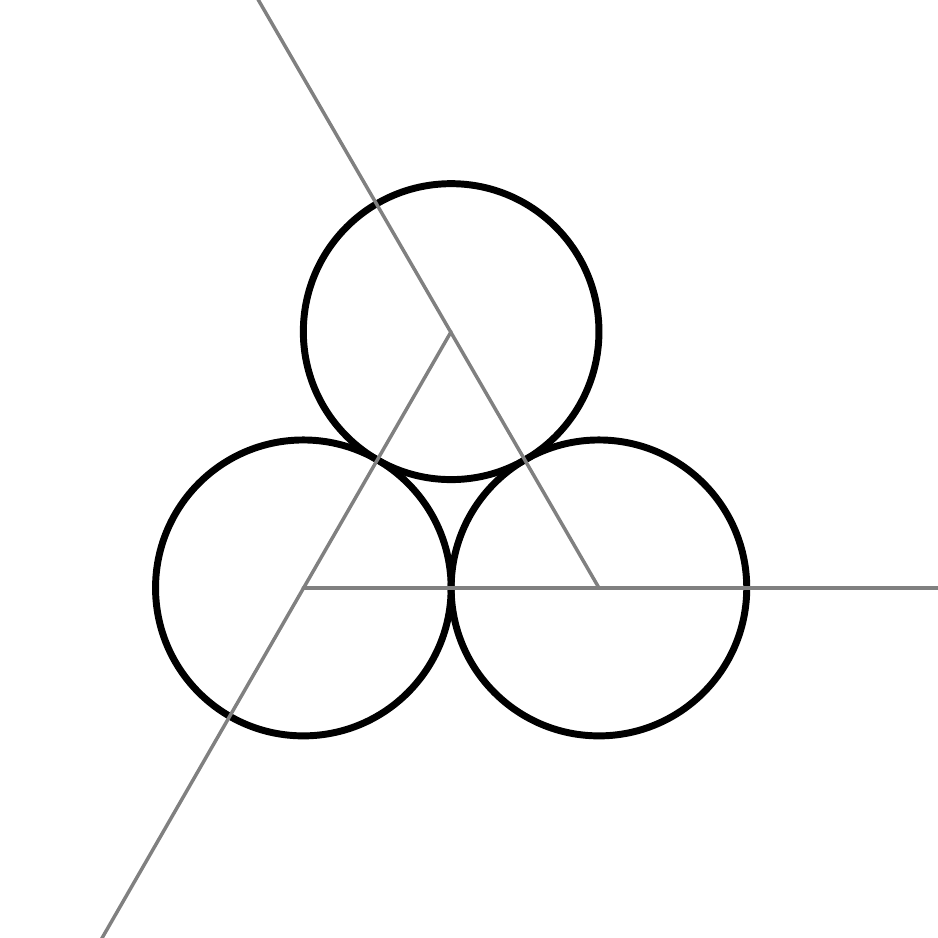}
\caption{Start with three circles of equal radii centered at the vertices of an equilateral triangle. Extend rays in order around the edges.}
\label{fig:construction1}
\end{subfigure}
\quad
\begin{subfigure}[b]{0.48\textwidth}
\includegraphics[width=\textwidth]{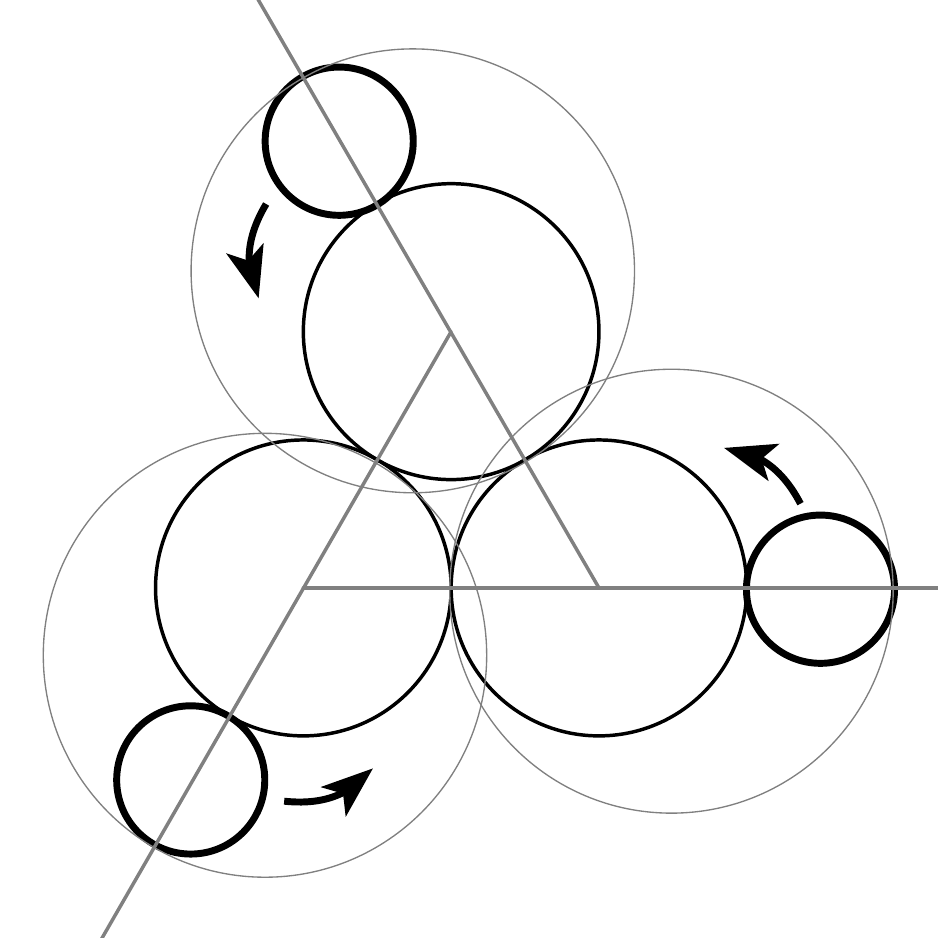}
\caption{Take three circles of equal radii at a fixed distance along each ray to form the outer face.}
\label{fig:construction2}
\end{subfigure}
\quad
\begin{subfigure}[b]{0.48\textwidth}
\includegraphics[width=\textwidth]{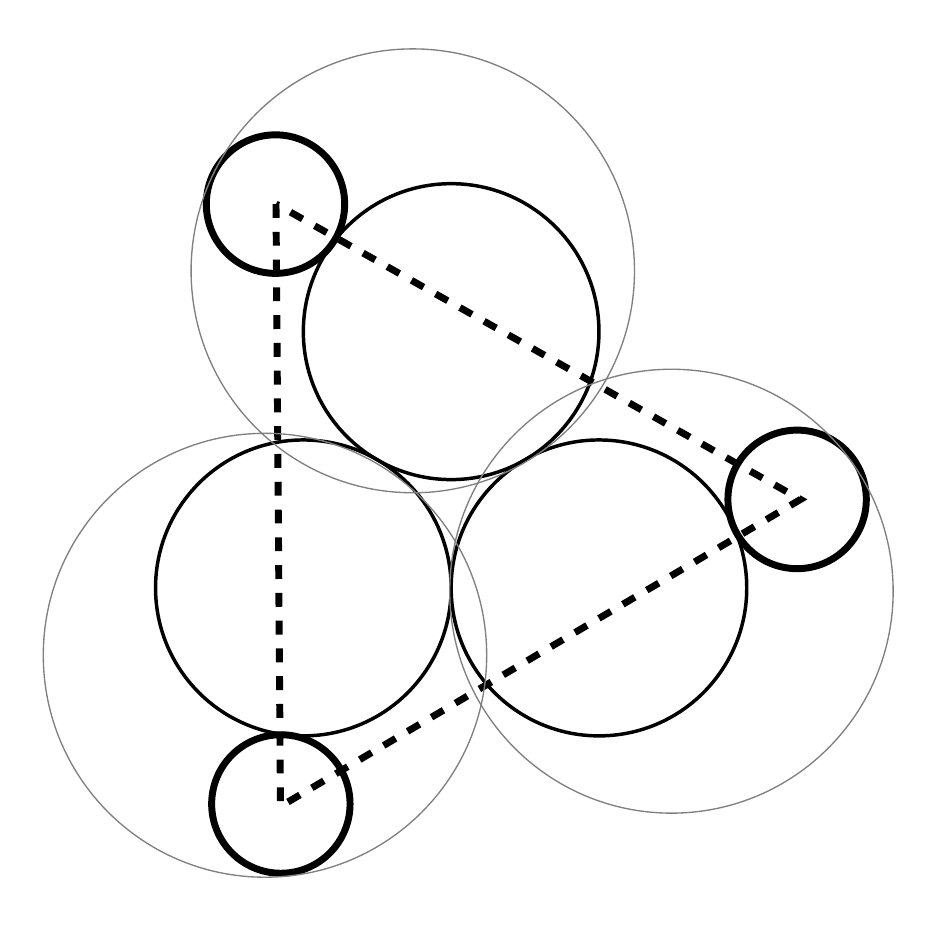}
\caption{Flow the outer circles symmetrically until the inversive distance reaches a minimum value. This is the critical circle realization.}
\label{fig:construction3}
\end{subfigure}
\quad
\begin{subfigure}[b]{0.48\textwidth}
\includegraphics[width=\textwidth]{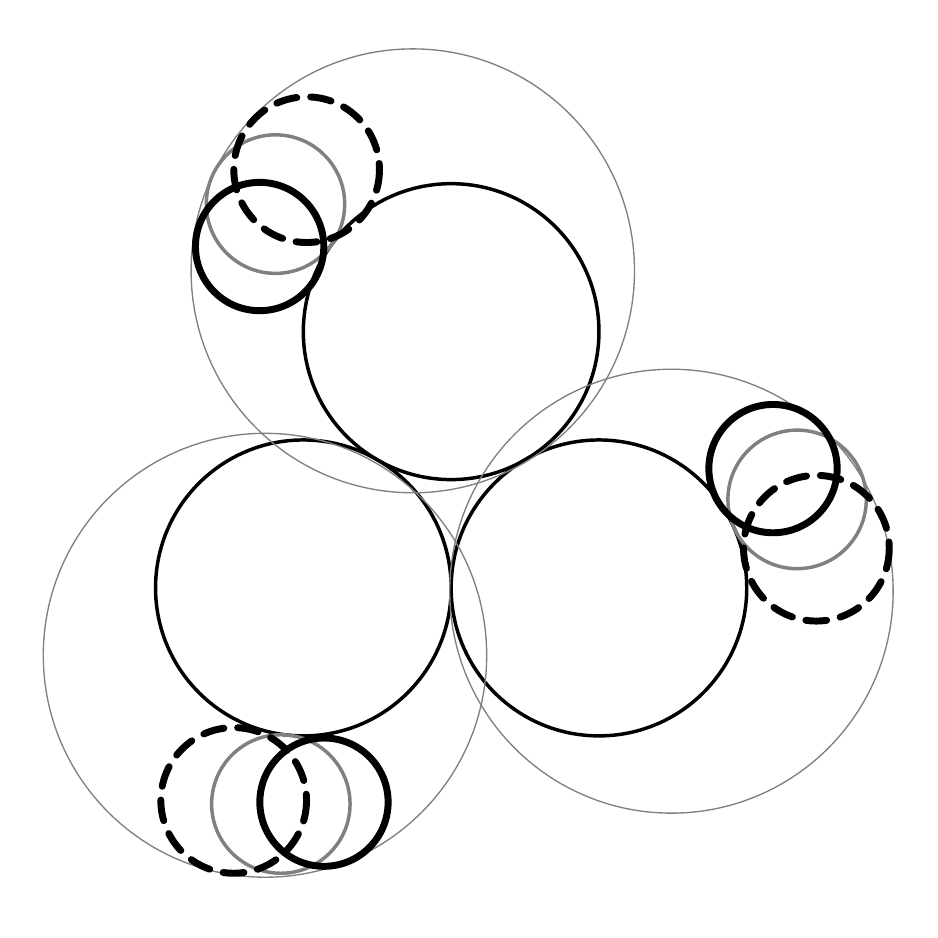}
\caption{Flow to find a pair of circle realizations (solid and dashed) near the critical packing that have equal inversive distances on the outer face. These are not M\"obius equivalent.}
\label{fig:construction4}
\end{subfigure}
\caption{The construction for a particular Ma-Schlenker pair.}\label{fig:construction}
\end{figure}

\subsection{The construction}\label{Section:Construction}
\begin{figure}
\begin{subfigure}[b]{0.3\textwidth}
\includegraphics[width=\textwidth]{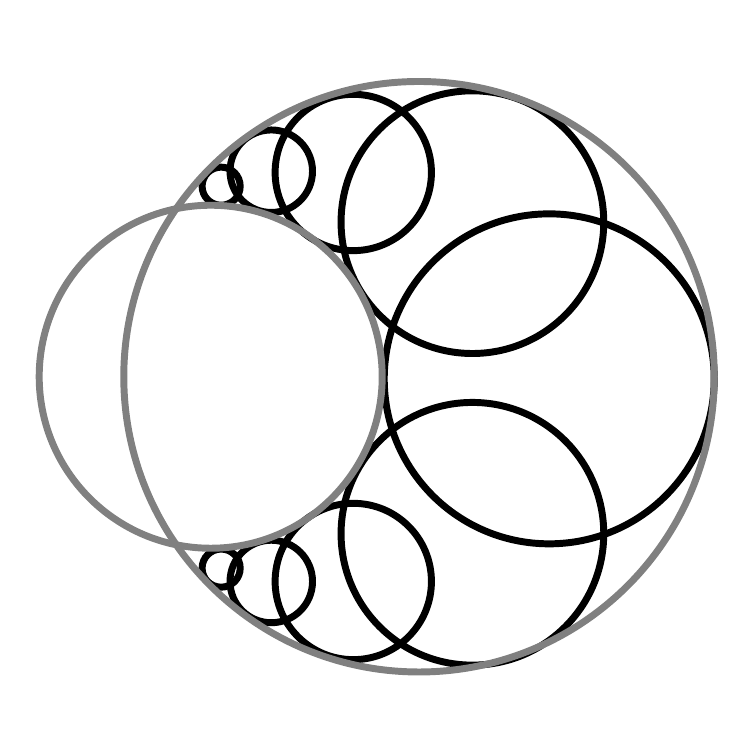}
\caption{Hyperbolic ($0 \leq a < 1$)}
\label{fig:flowhyperbolic}
\end{subfigure}
\quad
\begin{subfigure}[b]{0.3\textwidth}
\includegraphics[width=\textwidth]{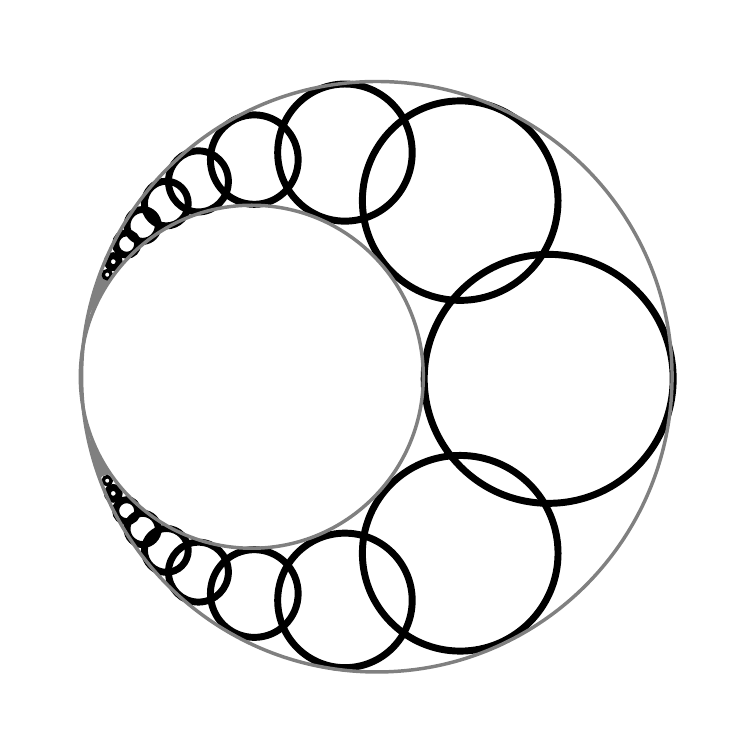}
\caption{Parabolic ($a = 1$)}
\label{fig:flowparabolic}
\end{subfigure}
\quad
\begin{subfigure}[b]{0.3\textwidth}
\includegraphics[width=\textwidth]{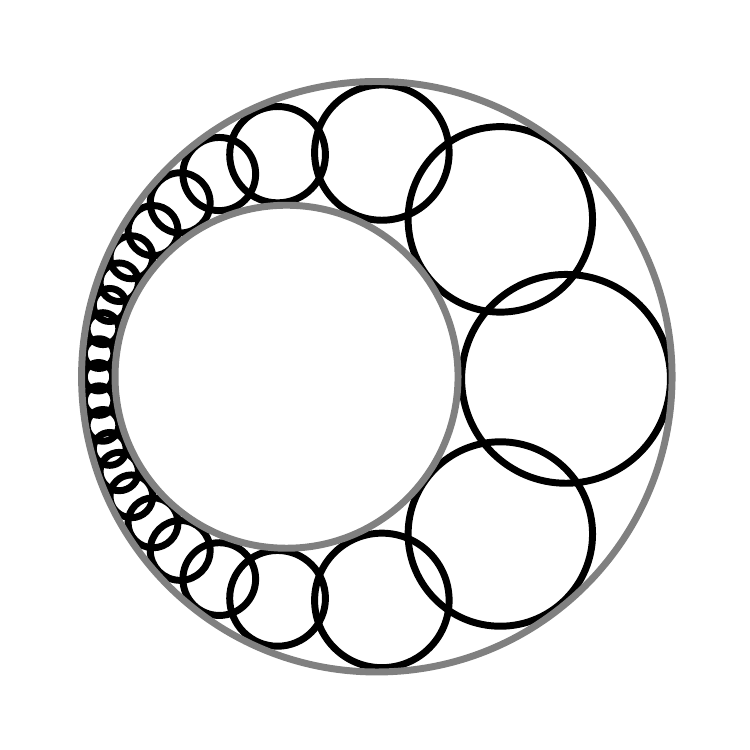}
\caption{Elliptic ($a > 1$)}
\label{fig:flowelliptic}
\end{subfigure}
\caption{M\"obius flows from the construction. The smaller gray circle is $A_1$ and the larger is $A_2$.}\label{fig:flow}
\end{figure}

We now describe a method for constructing Ma-Schlenker pairs. It begins with a construction of a $1$-parameter family of planar circle realizations of Ma-Schlenker octahedra $\mathscr{O}(a,b,c,d(t))$ where $a$, $b$, and $c$ are fixed while $d = d(t)$ varies. \Cref{fig:construction} presents a graphical description of this construction in a special case. Fix a non-negative number $a$ and choose three circles $C_{u}$, $C_{v}$, and $C_{w}$ in the complex plane that are centered at the respective vertices of an equilateral triangle $\Delta$ so that $\langle C_{u}, C_{v}\rangle = \langle C_{v}, C_{w}\rangle= \langle C_{w}, C_{u}\rangle =a$. Normalize by insisting that $\Delta$ have side length $2$ with circle $C_{u}$ centered on $z=-1$, $C_{v}$ centered on $z=1$, and $C_{w}$ centered on $z=\boldsymbol{i}\sqrt{3} $, placing the incenter of $\Delta$ at $z = \boldsymbol{i}/\sqrt{3}$. Let $\mathcal{A} = \mathcal{A}_{C_{u},C_{v}}$ be the coaxial family containing $C_{u}$ and $C_{v}$ and note that the line of centers of $\mathcal{A}$ is the real axis $\mathbb{R}$ and the radical axis of $\mathcal{A}$ is the imaginary axis $\mathbb{R}\boldsymbol{i}$. Choose an \textit{initial circle} $C$ centered on the positive real axis such that $c = \langle C_{v} , C \rangle \geq 0$ and let $b =\langle C_{u}, C\rangle$. Note that $b > c$. Let $\mu = \{\mu_{t}: t\in \mathbb{R}\}$ be the unit speed M\"obius flow determined by $\mathcal{A}$ that is counterclockwise on the circle $C_{v}$, and let $1 < x_{1} < x_{2}$ be the points of intersection of $C$ with the real axis. Let $A_{1}, A_{2} \in \mathcal{A}$ be the circles in the family $\mathcal{A}$ containing $x_{1}$ and $x_{2}$, respectively. For any real number $t$, $\mu_{t} (C)$ is a circle tangent to both $A_{1}$ and $A_{2}$ and sits between them, external to the disk bounded by $A_{1}$ and internal to the disk bounded by $A_{2}$. See \cref{fig:flow} for illustrations of the family $\{\mu_{t}(C) : t\in \mathbb{R}\}$ when the flow is hyperbolic ($0 \leq a <1$), parabolic ($ a = 1$), and elliptic ($a>1$). Since the flow $\mu$ preserves the individual circles of the coaxial system $\mathcal{A}$, $\langle C_{u}, \mu_{t} (C) \rangle = b$ and $\langle C_{v} , \mu_{t}(C) \rangle = c$ for all $t\in \mathbb{R}$.

Let $\mathfrak{r}$ be the counterclockwise rotation of the plane $\mathbb{C}$ through angle $2\pi/3$ with fixed point $\boldsymbol{i}/\sqrt{3}$, the incenter of the triangle $\Delta$. The circle configuration $\mathscr{C}(t)$ is defined as the collection
\begin{equation}
	\mathscr{C}(t) = \{ C_{u}, C_{v}, C_{w}, C_{w'}(t) = \mu_{t}(C), C_{u'}(t) = \mathfrak{r}(\mu_{t}(C)), C_{v'}(t) = \mathfrak{r}^{2}(\mu_{t}(C)) \}.
\end{equation}
The properties of $\mathscr{C}(t)$ of interest to us are, first, for all $t\in \mathbb{R}$ the configuration $\mathscr{C}(t)$ of circles in the plane is a circle realization for the Ma-Schlenker octahedron $\mathscr{O}(a,b,c,d(t))$, where 
\begin{equation}\label{EQ:d(t)}
	d(t) = \langle \mu_{t}(C), \mathfrak{r}(\mu_{t}(C)) \rangle,
\end{equation}
and second, for all $t\in \mathbb{R}$, $\mathscr{C}(t)$ has order three rotational symmetry about the incenter of $\Delta$ via the rotation $\mathfrak{r}$. Stereographic projection of $\mathscr{C}(t)$ to $\mathbb{S}^{2}$ defines a $1$-parameter family of circle realizations for $\mathscr{O}(a,b,c,d(t))$ in the $2$-sphere, though these are not \'a priori circle packings. These realizations possibly must be repositioned to a normalized position by applying M\"obius transformations of $\mathbb{S}^{2}$ to find appropriate intervals of $t$-values that produce the desired circle packings.
\begin{figure}
\includegraphics[width=0.4\textwidth]{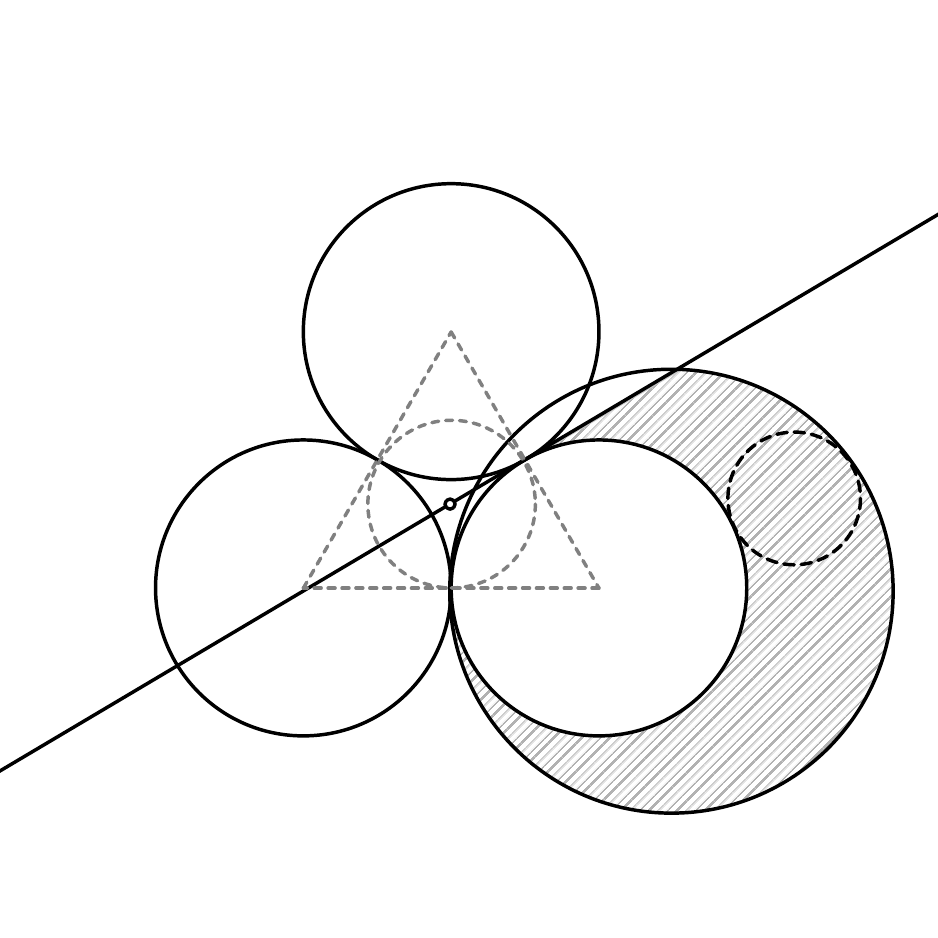}
\caption{The center of the critical circle must fall within the shaded region.}
\label{fig:criticalrange}
\end{figure}

\textit{Our goal is to find values of the inversive distance parameters $a > 1/2$ and $b > c \geq 0$ so that the function $d(t)$ has a critical value at a parameter value $t= \tau$, with $d(\tau)$ a local minimum for $d$, and for which the circle $\mu_{\tau}(C)$ is centered in the interior of the closed half-plane $\mathsf{H}$ containing $z = 1$, the center of $C_{v}$, and bordered by the line through $z=-1$, the center of $C_{u}$, and the incenter of $\Delta$; see \cref{fig:criticalrange}. Moreover, $\mu_{\tau}(C)$ should have positive inversive distance to the unique circle $O$ orthogonal to $C_{u}$, $C_{v}$, and $C_{w}$. The circle $O$ exists since $a >1/2$. The claim is that with these conditions satisfied, Ma-Schlenker pairs may be produced.}

 To see this, let $O'$ be the unique circle orthogonal to $C_{u'}(\tau)$, $C_{v'}(\tau)$, and $C_{w'}(\tau)$, which exists since $d(\tau) > 1/2$.\footnote{By an application of Equation~\ref{EQ:theta}.} By the order three rotational symmetry of $\mathscr{C}(t)$, $O$ and $O'$ are concentric, both centered at $\boldsymbol{i}/\sqrt{3}$. Let $\widetilde{\mathscr{C}}(t)$ be the image of $\mathscr{C}(t)$ under stereographic projection to the $2$-sphere followed by a M\"obius transformation so that (1) the respective images of the circles $O$ and $O'$ are latitudinal circles, the first centered on the south pole and contained in the southern hemisphere, and the second on the north pole and contained in the northern hemisphere, and (2) the circles $L$ and $L'$ have the same radius, $L$ centered on the south pole and $L'$ on the north, where $L$ is the latitudinal circle containing the centers of the projections of $C_{u}$, $C_{v}$, and $C_{w}$, and $L'$ is the latitudinal circle containing the centers of the projections of $C_{u'}(\tau)$, $C_{v'}(\tau)$, and $C_{w'}(\tau)$. Note that the order three rotational symmetry of $\mathscr{C}(t)$ translates to an order three rotational symmetry of $\widetilde{\mathscr{C}}(t)$ by a rotation of $\mathbb{S}^{2}$ about the axis through the north and south poles. 

\begin{figure}
\begin{subfigure}[b]{0.48\textwidth}
\includegraphics[width=\textwidth]{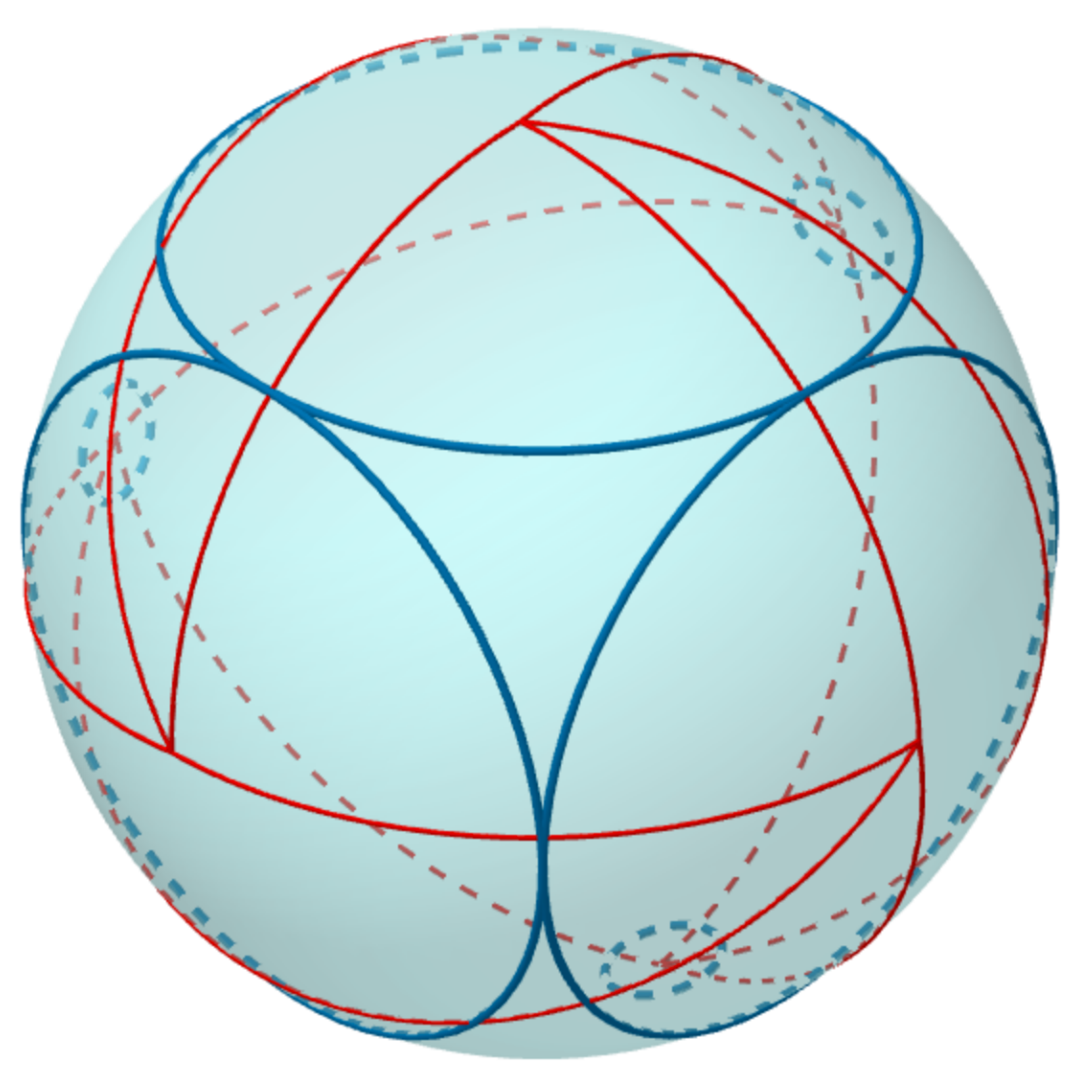}
\caption{A view from the south pole.}
\label{fig:3dSouthPole}
\end{subfigure}
\quad
\begin{subfigure}[b]{0.48\textwidth}
\includegraphics[width=\textwidth]{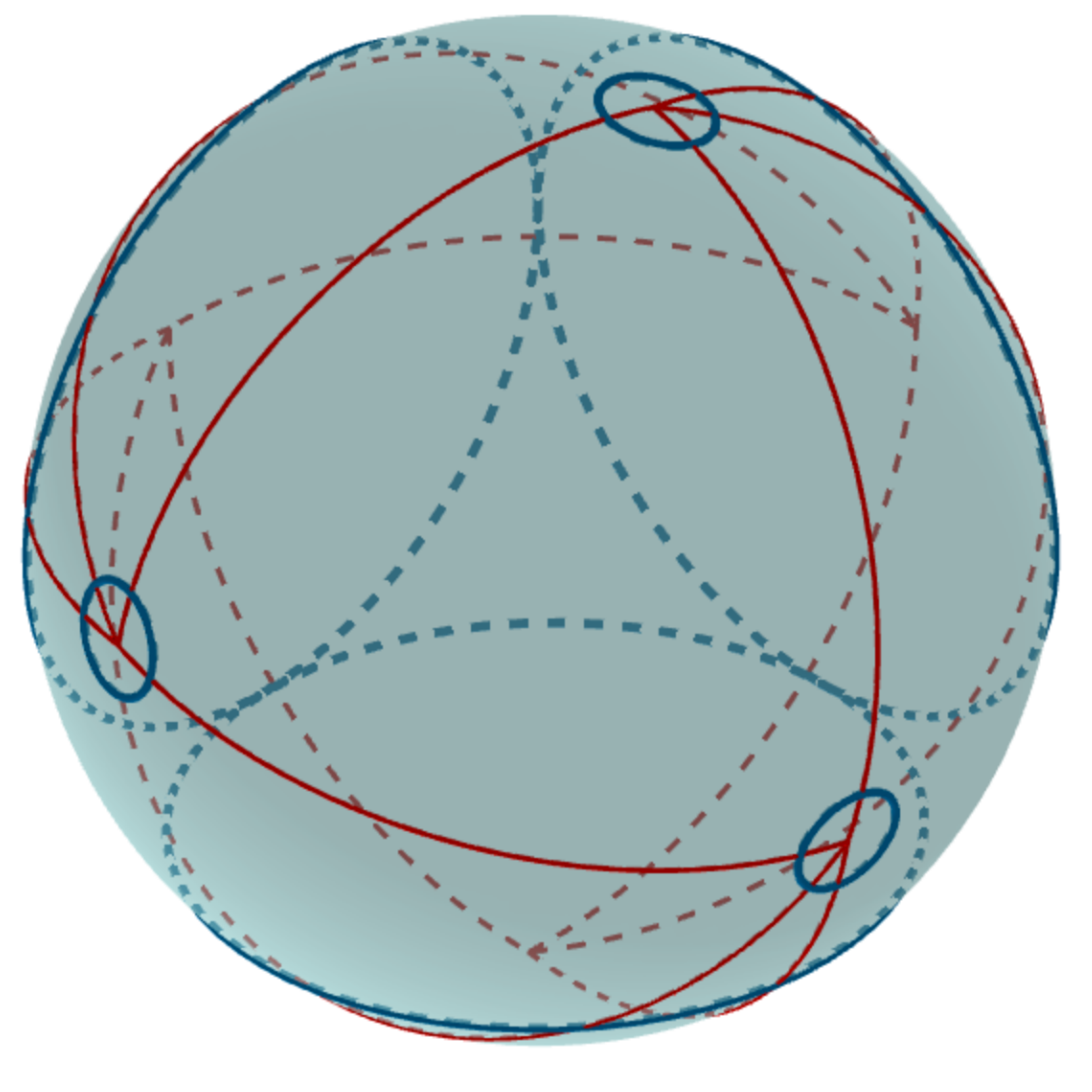}
\caption{A view from the north pole.}
\label{fig:3dNorthPole}
\end{subfigure}
\quad
\begin{subfigure}[b]{0.48\textwidth}
\includegraphics[width=\textwidth]{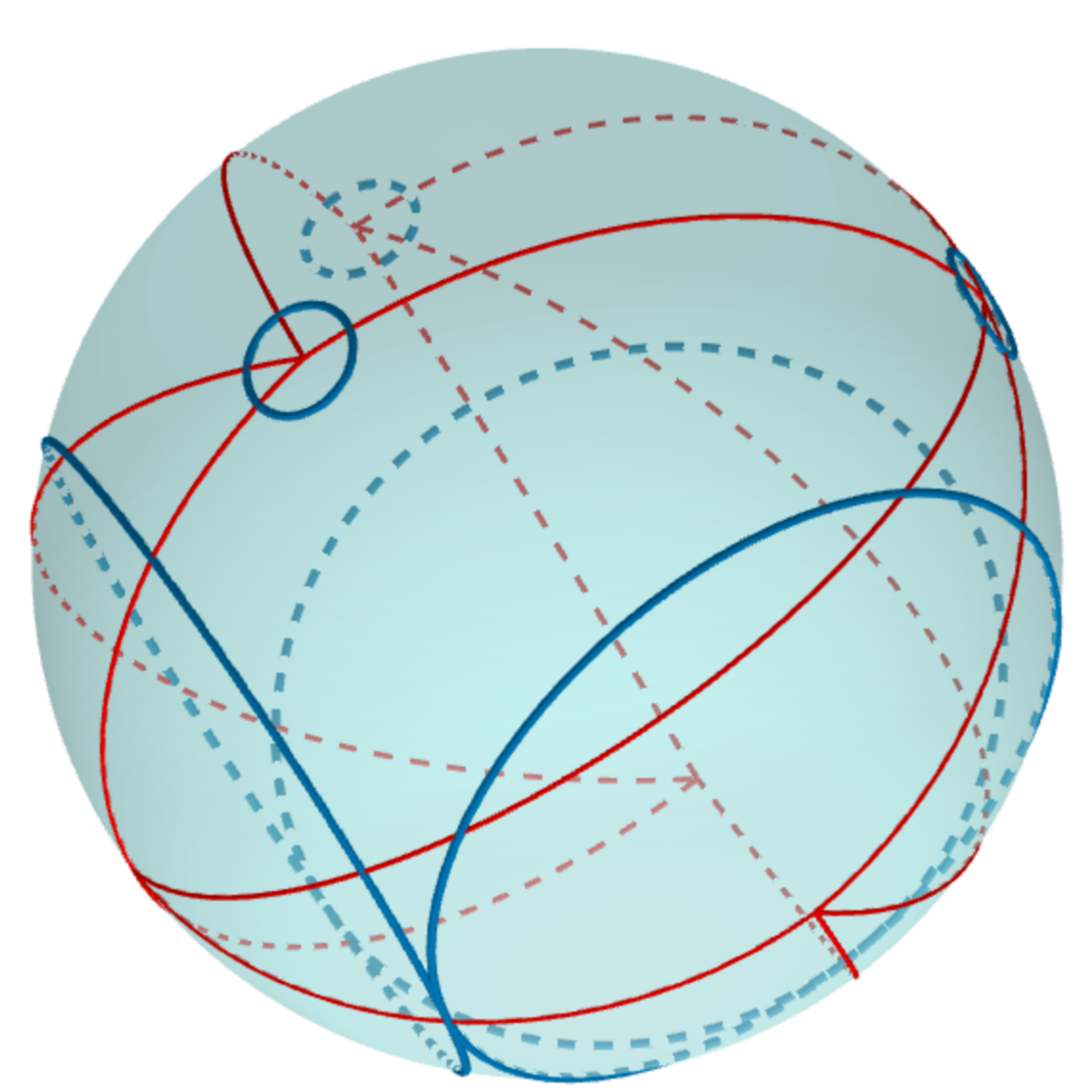}
\caption{A view from slightly north of the equator.}
\label{fig:3dView1}
\end{subfigure}
\quad
\begin{subfigure}[b]{0.48\textwidth}
\includegraphics[width=\textwidth]{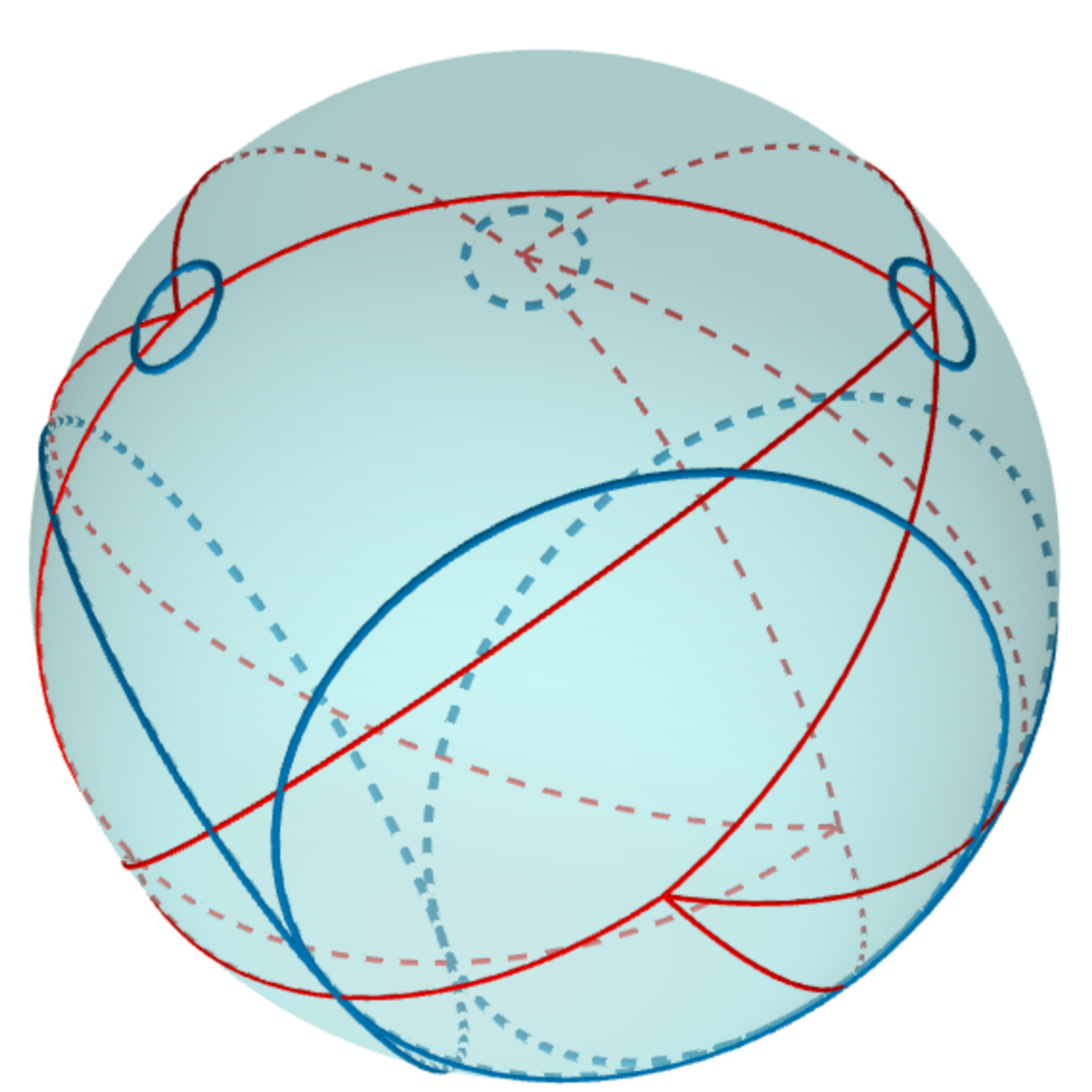}
\caption{A slight rotation from \cref{fig:3dView1}.}
\label{fig:3dView2}
\end{subfigure}
\caption{A 3D visualization of a critical Ma-Schlenker octahedron showing the triangulation.}\label{fig:3d}
\end{figure}

Our claim is that there is an open interval $J$ of $t$-values containing $\tau$ for which each $\widetilde{\mathscr{C}}(t)$, for $t\in J$, is a circle packing. Since then the Ma-Schlenker \textit{c}-octahedron $\widetilde{\mathscr{C}}(\tau)$ provides a minimum value for $d(t)$, variation of $t$ about $\tau$ produces (possibly with further restrictions on $a$, $b$, and $c$) pairs $t < \tau$ and $t' > \tau$ in $J$ for which $\widetilde{\mathscr{C}}(t)$ and $\widetilde{\mathscr{C}}(t')$ form a Ma-Schlenker pair, two circle packings for $\mathscr{O}(a,b,c,d)$, where $d= d(t) = d(t')$, that fail to be inversive equivalent. 
\begin{Definition}
The circle realization $\mathscr{C}(\tau)$ is said to be \textit{critical}, the circle packing $\widetilde{\mathscr{C}}(\tau)$ is a \textit{critical Ma-Schlenker \textit{c}-octahedron}, and the corresponding edge-labeled triangulation $\mathscr{O}(a,b,c,d(\tau))$ is a \textit{critical Ma-Schlenker octahedron}. See \cref{fig:3d} for a 3D visualization of a critical Ma-Schlenker \textit{c}-octahedron.
\end{Definition}

The claim of the preceding paragraph, that there is an open interval $J$ of $t$-values containing $\tau$ for which each $\widetilde{\mathscr{C}}(t)$, for $t\in J$, is a circle packing follows from the fact that the circle $\mu_{\tau}(C)$ is centered in the interior of the closed half-plane $\mathsf{H}$. The point is that three points $u$, $v$, and $w$ equally spaced on $L$ and $u'$, $v'$ and $w'$ equally spaced on $L'$ with $u$ and $w'$ on the same meridian, $v$ and $u'$ on the same meridian, and $w$ and $v'$ on the same meridian, cut out an octahedral triangulation of $\mathbb{S}^{2}$ when points are connected according to the Ma-Schlenker octahedral pattern. Assuming the points are ordered counterclockwise when viewed from the north pole, if now the points along $L'$ are rotated counterclockwise by an  angle strictly between $0$ and $\pi$, then this triangulation ``stretches'' to form a triangulation with vertices $u$, $v$, and $w$, and the rotated $u'$, $v'$, and $w'$. It is only when the rotation reaches $\pi$ radians that there is ambiguity as then $u$ and $w'$ are antipodal, and then as $w'$ is rotated past $\pi$ radians, the geodesic arc connecting $u$ to $w'$ moves to the other side of the north pole, and the triangle $uvw'$ now contains the north pole, as do $vwu'$ and $wuv'$. At this point we have lost the triangulation. The point of the half-plane is that the center of $\mu_{\tau}(C)$ in the interior of $\mathsf{H}$ guarantees that the position of $w'$ is obtained in this manner by a positive counterclockwise rotation strictly less than $\pi$.

We have used the restriction $a > 1/2$ to ensure the existence of a circle mutually orthogonal to the three circles $C_{u}$, $C_{v}$, and $C_{w}$. This in turn is used to ensure that stereographic projection followed by an appropriate M\"obius transformation produces, not just a realization, but a circle packing for a Ma-Schlenker octahedron. This is done for cenvenience of argument and we mention that even if $a \leq 1/2$, when there is no mutually orthogonal circle to $C_{u}$, $C_{v}$, and $C_{w}$, Ma-Schlenker circle packings are still possible.

\section{Ma-Schlenker \textit{c}-Octahedra---the Devil of the Details}\label{Section:Details}
\begin{figure}
\begin{subfigure}[t]{0.3\textwidth}
\includegraphics[width=\textwidth]{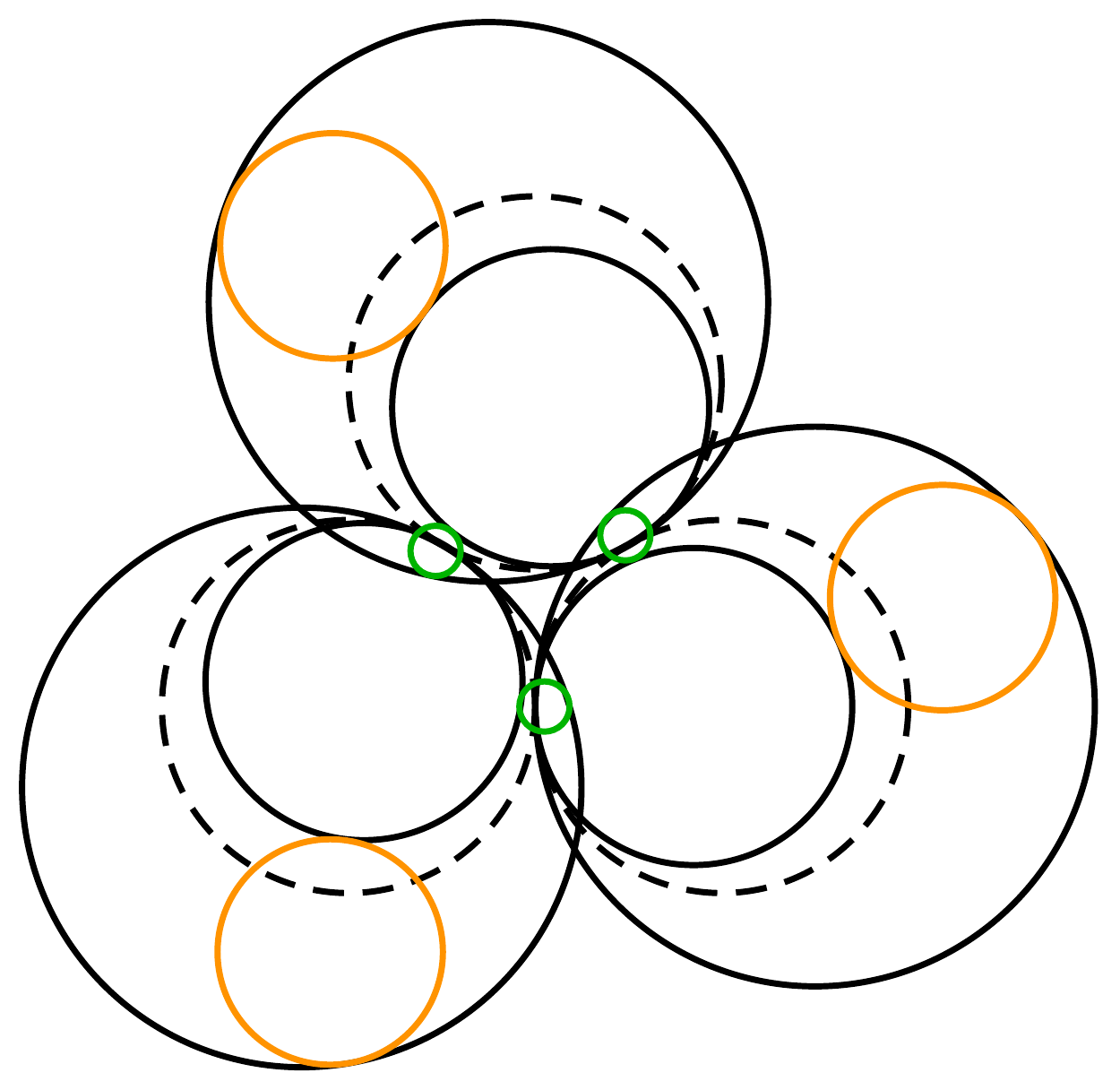}
\caption{Parabolic flow.}
\label{fig:flowgraph_para1}
\end{subfigure}
\quad
\begin{subfigure}[t]{0.6\textwidth}
\includegraphics[width=\textwidth]{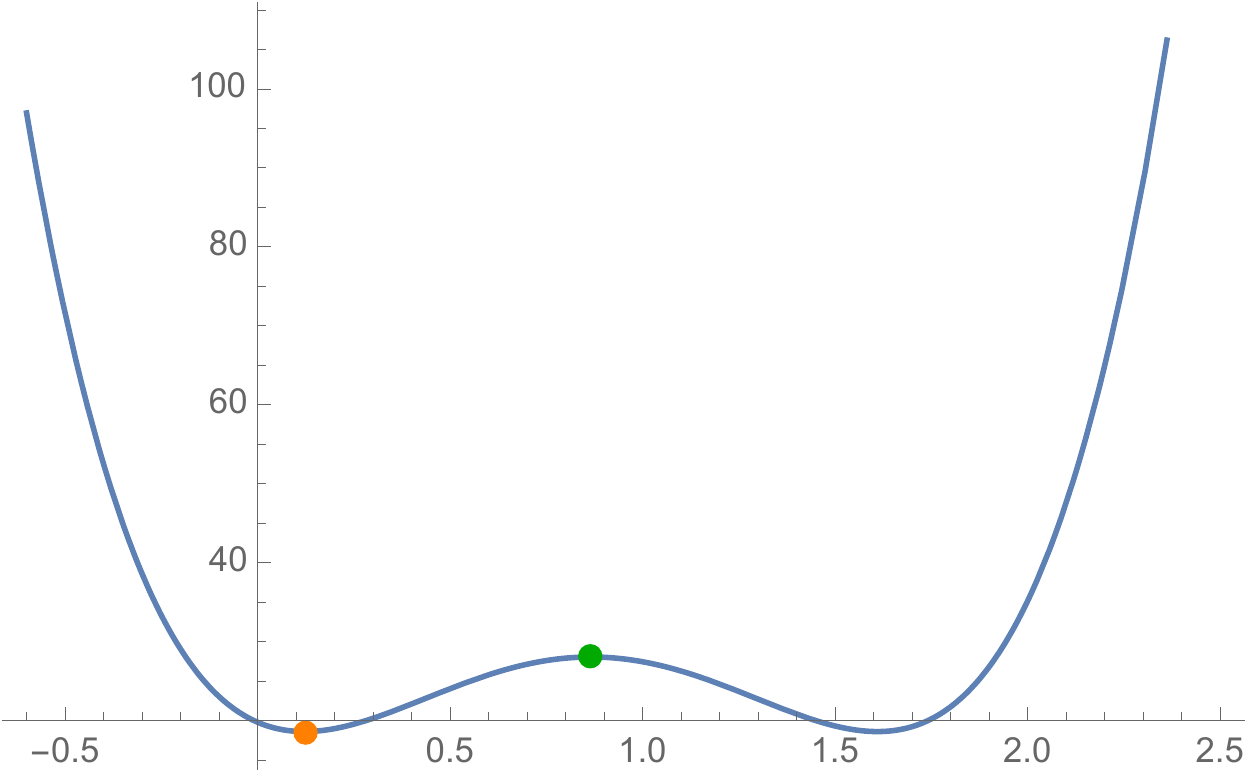}
\caption{Graph of the inversive distance function $d(t)$.}
\label{fig:flowgraph_para2}
\end{subfigure}
\caption{A parabolic flow. The parameters are $x_1 = 1.7$, $x_2 = 3$. The highlighted minimum and (local) maximum occur at $\tau = t_{min} = 0.121766$ (orange) and $m = t_{max} = 0.866025$ (green), and the inversive distances are $d(t_{min}) = 18.6065$ and $d(t_{max}) = 28.051$. (Computed values are rounded.)}
\label{fig:flowgraph_para}
\end{figure}

\begin{figure}
\begin{subfigure}[t]{0.35\textwidth}
\includegraphics[width=\textwidth]{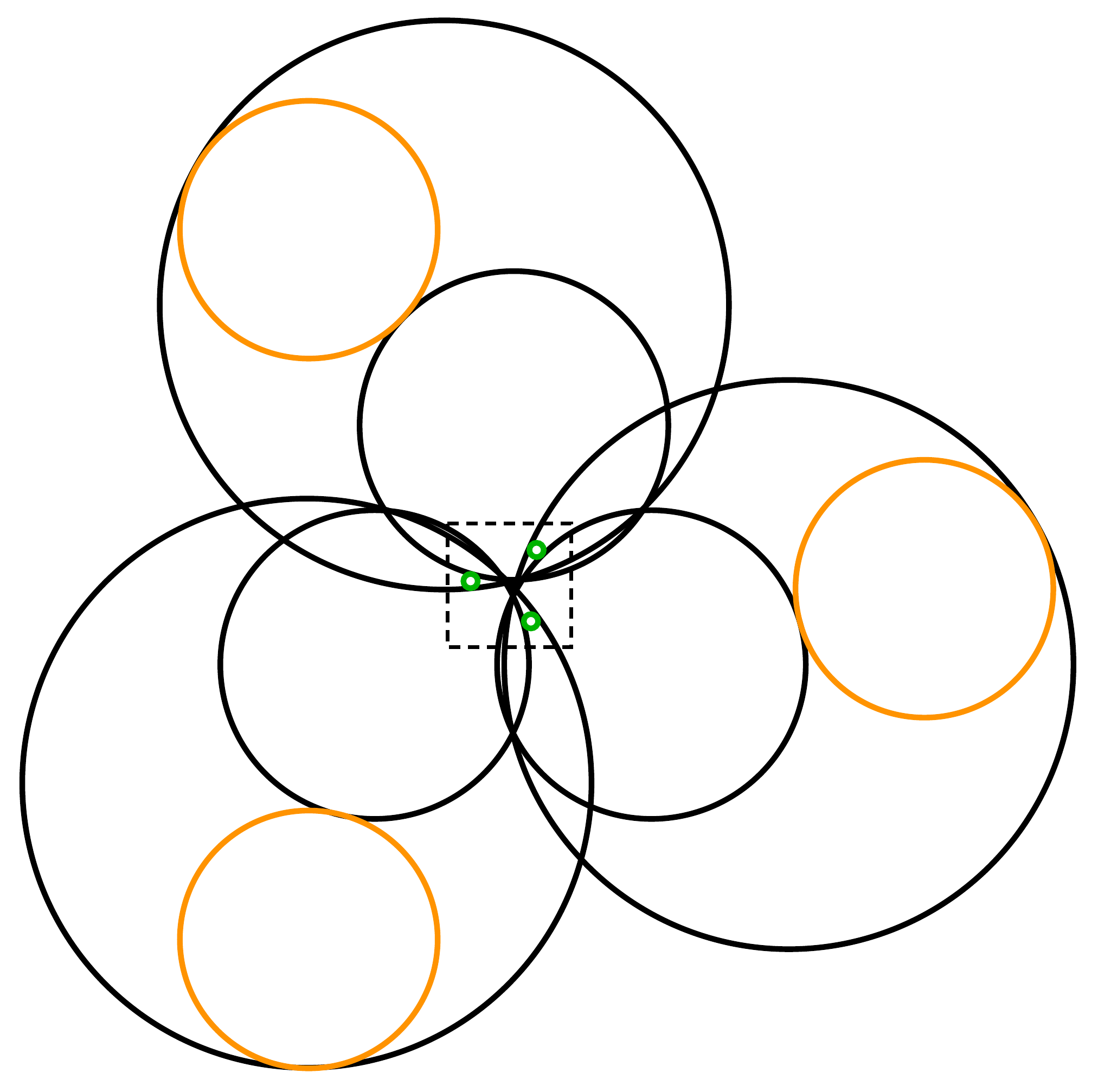}
\caption{Hyperbolic flow.}
\label{fig:flowgraph_hyp1}
\end{subfigure}
\quad
\begin{subfigure}[t]{0.2\textwidth}
\includegraphics[width=\textwidth]{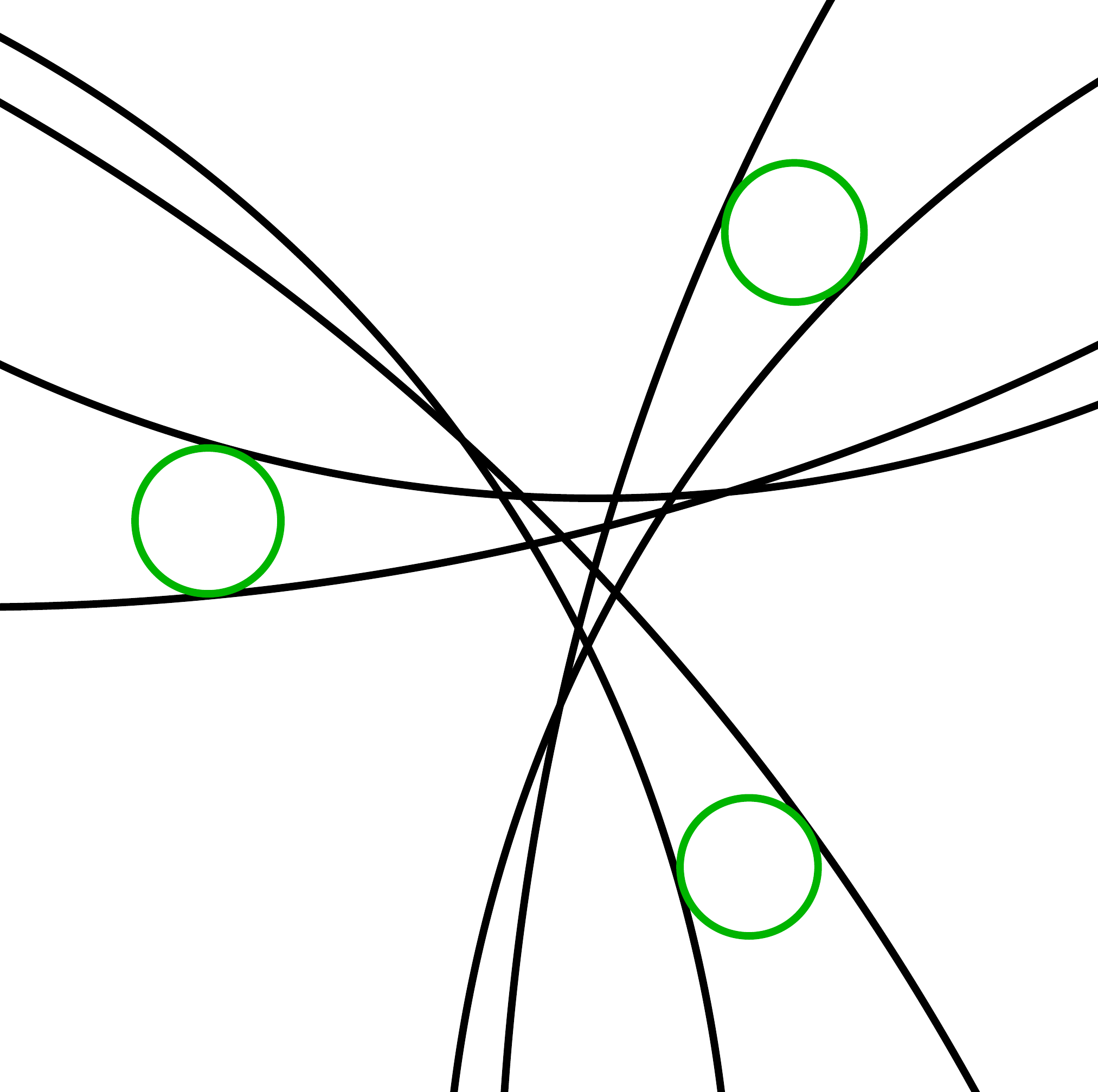}
\caption{Detail of the dotted box.}
\label{fig:flowgraph_hyp2}
\end{subfigure}
\quad
\begin{subfigure}[t]{0.35\textwidth}
\includegraphics[width=\textwidth]{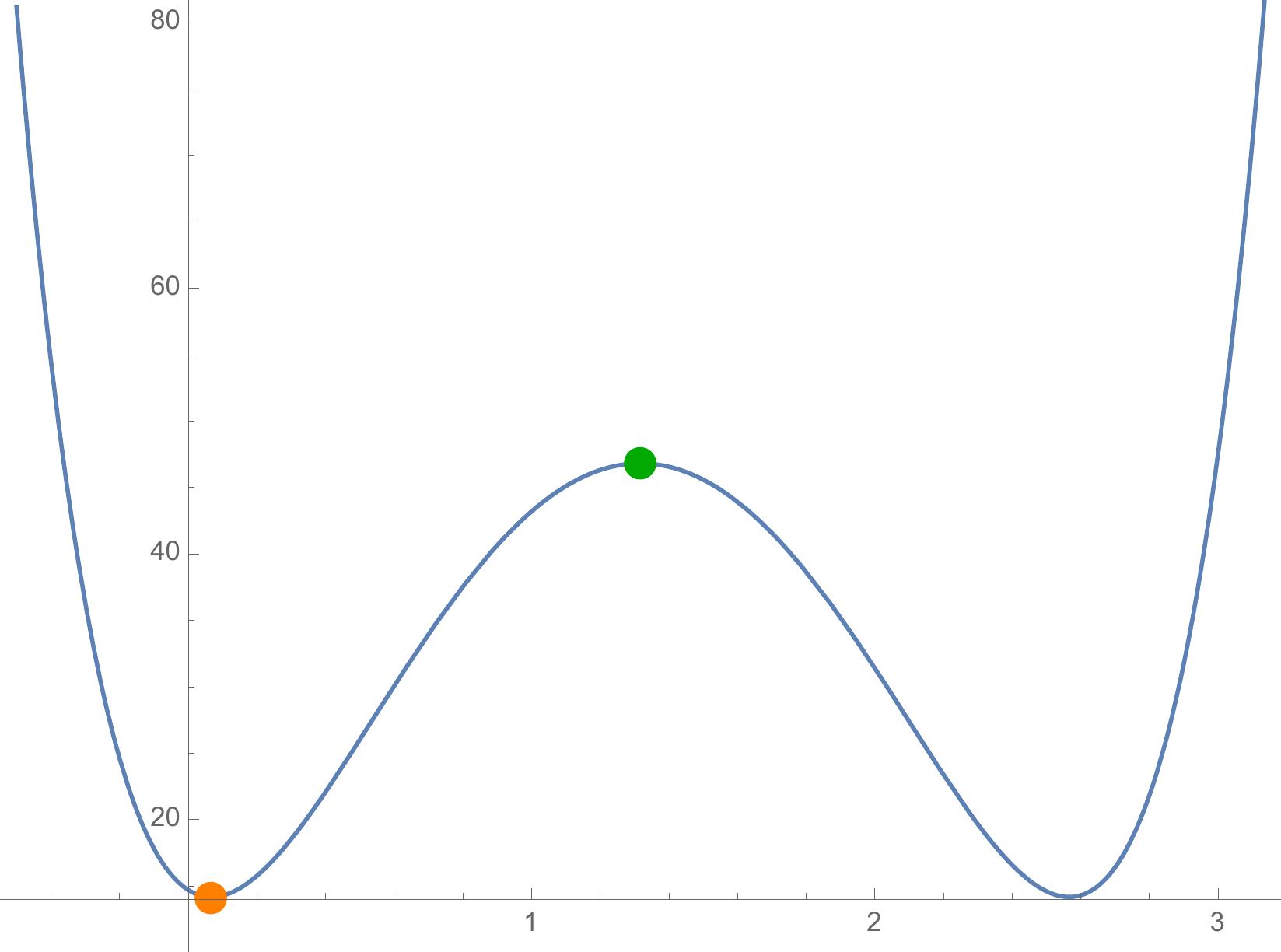}
\caption{Graph of the inversive distance.}
\label{fig:flowgraph_hyp3}
\end{subfigure}
\caption{Hyperbolic flow. The parameters are $x_1 = 2.11803$, $x_2 = 4.06155$, and $y = 0.5$. The highlighted minimum (orange) occurs at $\tau = t_{min} = 0.06782$ and the local maximum (green) at $m = t_{max} = 1.31696$. The inversive distances are $d(t_{min}) = 14.1647$ and $d(t_{max}) = 46.8136$.}
\label{fig:flowgraph_hyp}
\end{figure}

\cref{fig:flowgraph_para,fig:flowgraph_ellip,fig:flowgraph_hyp} show examples of the graphs of $d(t)$ for various values of $a$ and choices of initial circle $C$. The Ma-Schlenker pairs arise from pairs $t$ and $t'$ near $\tau$ and at the same horizontal level on the graph of $d(t)$. A precise description of the derivation of these graphs is given subsequently, after we make some observations.
%
%

\subsection{The parabolic case}
Here $a=1$, the circles $C_{u}$, $C_{v}$ and $C_{w}$ are mutually tangent, and the flow is parabolic. An example is shown in \cref{fig:flowgraph_para1}, where the initial circle $C$ meets the real axis at $x_{1} = 1.7$ and $x_{2} = 3$, with $b = \langle C_{u}, C\rangle = 7.538$ and $c = \langle C_{v}, C\rangle = 0.308$. The three circles $C_{u}$, $C_{v}$ and $C_{w}$ are the dotted circles. The critical circle $C(\tau)$, in orange, gives a minimum value for $d(t)$ and is centered in the half-plane $\mathsf{H}$. The graph of $d(t)$ is shown in \cref{fig:flowgraph_para2}, whose shape is typical of all the examples with $a=1$. The general characteristics of the parabolic flow are the same as for the hyperbolic, which is described in detail in the next paragraph.

\subsection{The hyperbolic case}
Here $1/2 < a\leq 1$ and the pairs of circles among $C_{u}$, $C_{v}$, and $C_{w}$ meet at angle $\pi/3 \geq \theta = \cos^{-1} a \geq 0$. The center of the unique circle $O$ orthogonal to $C_{u}$, $C_{v}$ and $C_{w}$ is the incenter $\boldsymbol{i}/\sqrt{3}$ of $\Delta$, and lies in the bounded interstice formed by $C_{u}$, $C_{v}$ and $C_{w}$. As $t\to \pm\infty$, the inversive distance parameter $d(t) \to + \infty$. The absolute minimum value of $d(t)$ occurs at two distinct $t$-values $0 < \tau < \tau'$, which implies that both circles, $C(\tau) =\mu_{\tau}(C)$ and $C(\tau') = \mu_{\tau'}(C)$, are centered in the first quadrant of the complex plane. Between these lies a parameter value $m$ at which $d(t)$ obtains a local maximum. The circle $C(m) = \mu_{m}(C)$ is orthogonal to $O$ and inversion $I_{O}$ through the circle $O$ preserves $C(m)$ and exchanges circles $C(\tau)$ and $C(\tau')$. The center of $C(\tau)$ lies outside of $O$, and that of $C(\tau')$ inside. More generally, the inversion $I_{O}$ generates a symmetry of the graph in the following way. Since $O$ is in the orthogonal complement $\mathcal{A}_{C_{u},C_{v}}^{\perp}$, the inversion $I_{O}$ preserves the family $\mathcal{A} =\mathcal{A}_{C_{u},C_{v}}$, merely inverting each circle of $\mathcal{A}$ to itself. In particular, since $A_{1}$ and $A_{2}$ are members of $\mathcal{A}$ and form the envelope of the family $\mathcal{B} = \{\mu_{t}(C) : t \in \mathbb{R}\}$ of circles generated by the M\"obius flow $\mu$, $I_{O}(A_{1}) = A_{1}$ and $I_{O}(A_{2}) = A_{2}$ form the envelope of the circles of the family $I_{O}(\mathcal{B})$, implying that $I_{O}$ preserves the family $\mathcal{B}$. It follows that $I_{O}$ generates an involution of $\mathbb{R}$. Indeed, for each $t\in \mathbb{R}$, let $t'$ be the real value for which $\mu_{t'}(C) = I_{O}(\mu_{t}(C))$. Then $t'' =t$ with fixed point $m = m'$. This generates a symmetry of the graph of $d(t)$ with $d(t) = d(t')$ since
\begin{align*}
	 d(t) = \langle \mu_{t}(C), \mathfrak{r}(\mu_{t}(C)) \rangle &=  \langle I_{O}(\mu_{t}(C)), I_{O}(\mathfrak{r}(\mu_{t}(C)) )\rangle \\ &= \langle I_{O}(\mu_{t}(C)), \mathfrak{r}(I_{O}(\mu_{t}(C))) \rangle = \langle \mu_{t'}(C), \mathfrak{r}(\mu_{t'}(C)) \rangle =d(t'),
\end{align*}
since the inversion $I_{O}$ commutes with the rotation $\mathfrak{r}$, as both are centered at $\boldsymbol{i}/\sqrt{3}$.

\cref{fig:flowgraph_hyp} presents an example where $a= 0.6$ with the angle $\theta$ of intersection of $C_{u}$ and $C_{v}$ equal to $2 \tan^{-1} (1/2)\cong 0.9273$. The three circles $C_{u}$, $C_{v}$ and $C_{w}$ are suppressed in the figure, and the respective values of $b$ and $c$ are $b=6.689$ and $c=1$. The circles $C_{u}$ and $C_{v}$ meet at $\pm y \boldsymbol{i}$ where $y = 0.5$. In the close up of \cref{fig:flowgraph_hyp2}, the three local maximum circles (green) and the six enveloping circles are orthogonal to $O$ (not shown).

\begin{figure}
\begin{subfigure}[t]{0.3\textwidth}
\includegraphics[width=\textwidth]{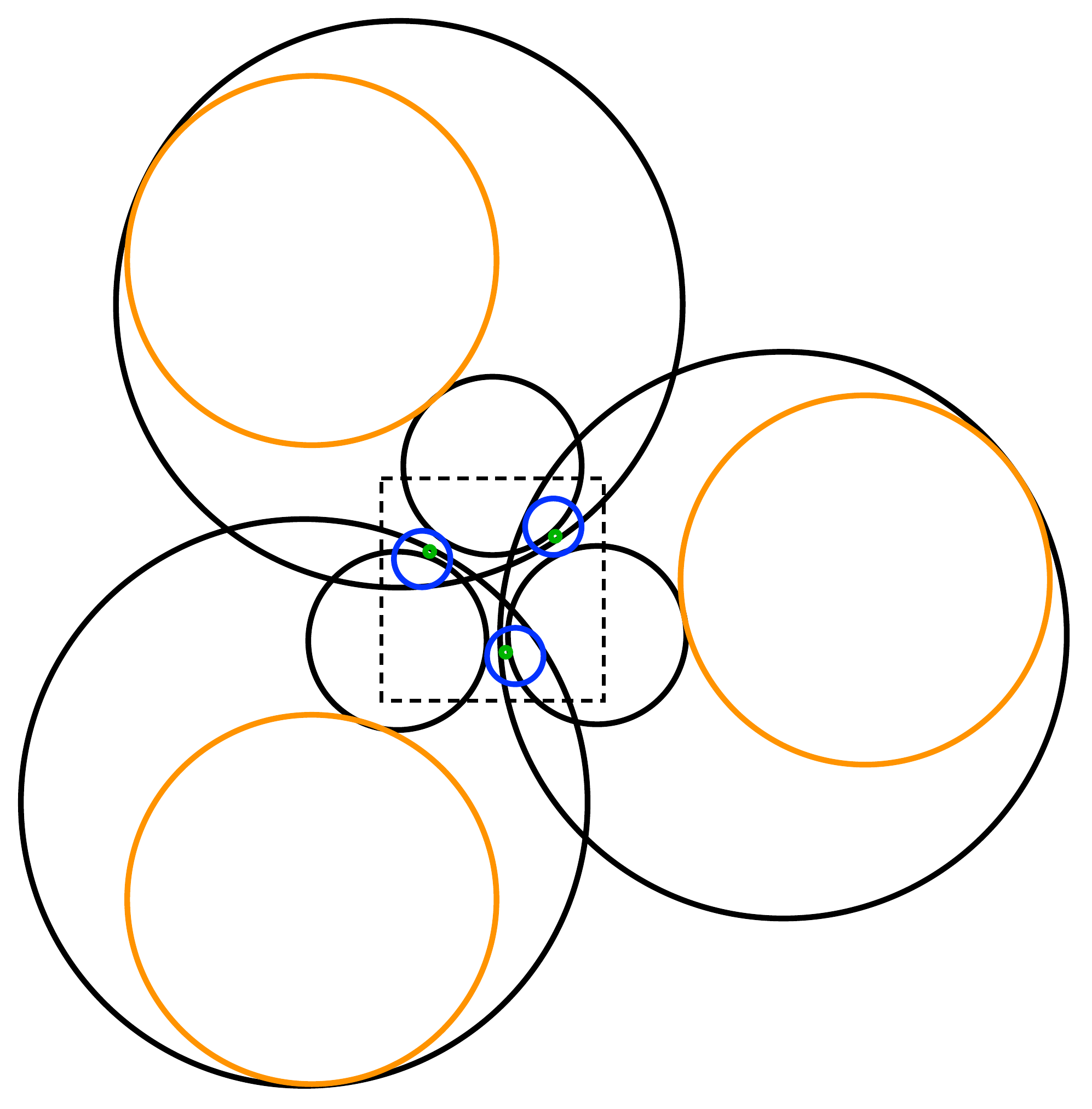}
\caption{Elliptic flow.}
\label{fig:flowgraph_ellip1}
\end{subfigure}
\quad
\begin{subfigure}[t]{0.2\textwidth}
\includegraphics[width=\textwidth]{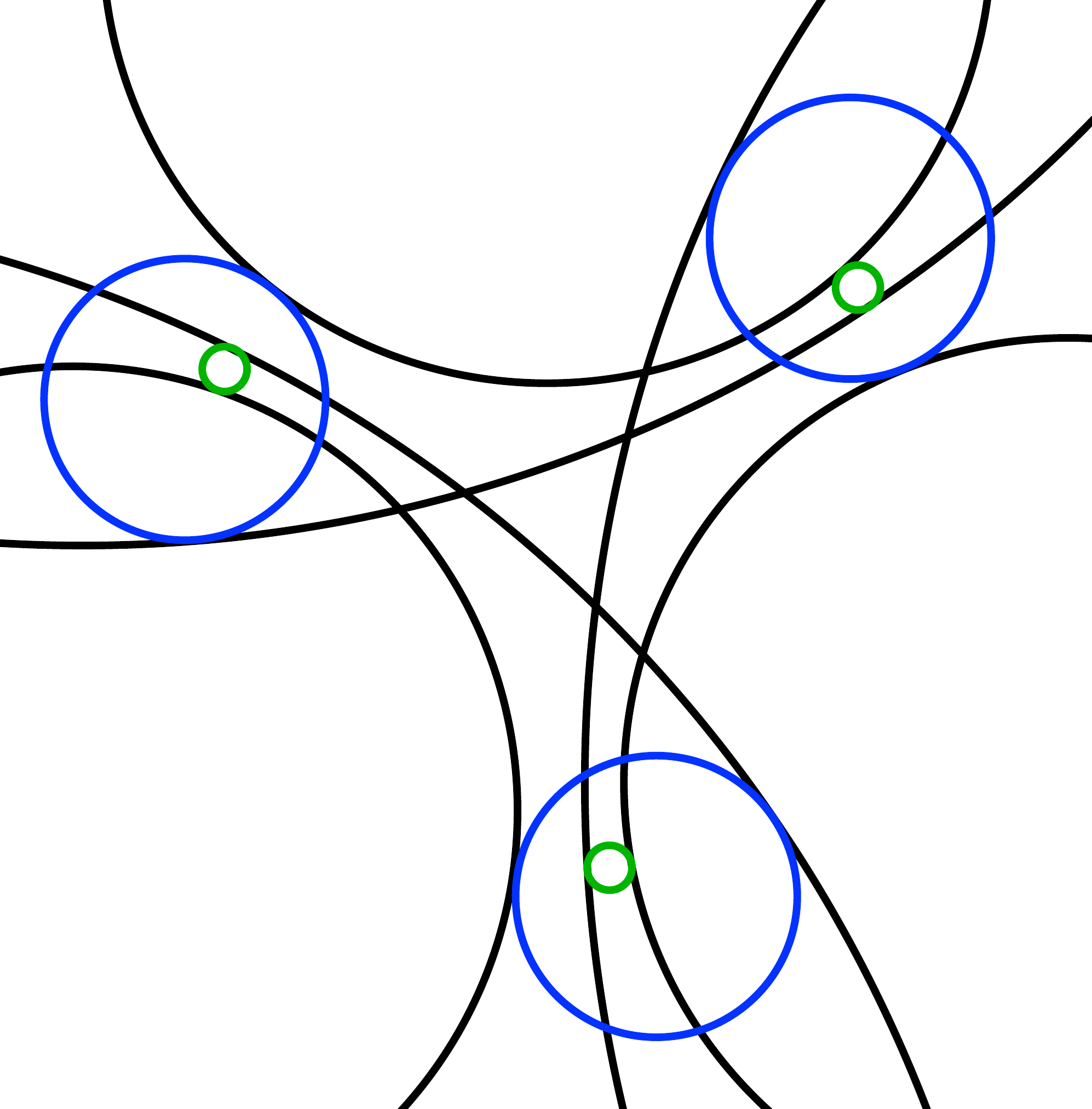}
\caption{Detail of the dotted box.}
\label{fig:flowgraph_ellip2}
\end{subfigure}
\quad
\begin{subfigure}[t]{0.4\textwidth}
\includegraphics[width=\textwidth]{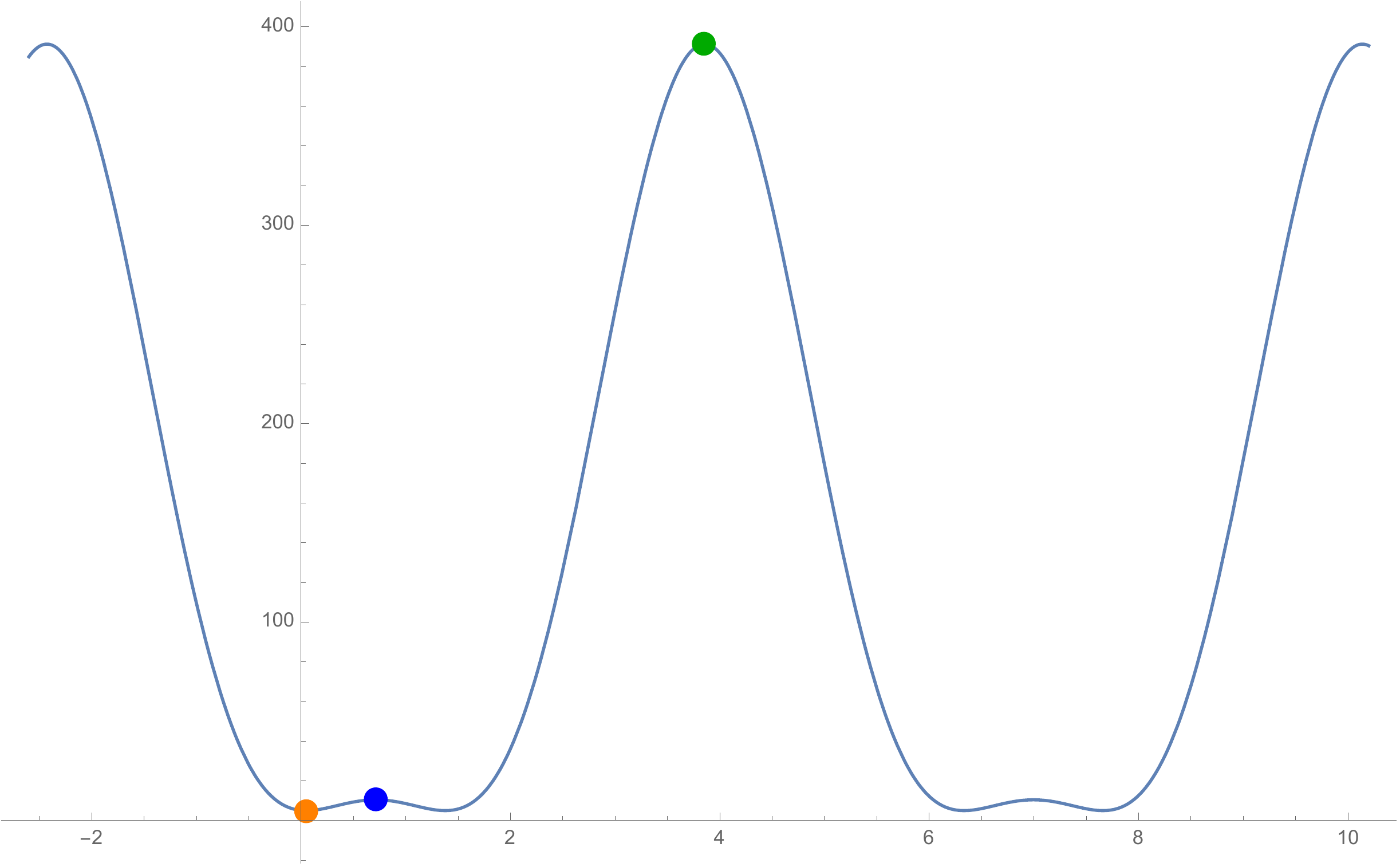}
\caption{Graph of the inversive distance, which is periodic of period $2\pi$.}
\label{fig:flowgraph_ellip3}
\end{subfigure}
\caption{An elliptic flow. The parameters are $x_1 = 2$, $x_2 = 6$, and $y = 0.5$. The highlighted minimum (orange) occurs at $\tau = t_{min} = 0.0506$, the local maximum (blue) at $m = t_{max}^{(1)} = 0.7137$, and the absolute maximum (green) at $M = t_{max}^{(2)} = 3.85532$. The inversive distances are $d(t_{min}) = 5.000$, $d(t_{max}^{(1)}) = 10.4245$, and $d(t_{max}^{(2)}) = 391.247$.}
\label{fig:flowgraph_ellip}
\end{figure}
\subsection{The elliptic case}
Here $a > 1$ and the facts are similar to those of the hyperbolic case. The primary difference is that $d(t)$ now is periodic of period $\omega = 2\pi / |\lambda|$, where $|\lambda|$ is the speed of the elliptic flow. There are parameter values $0 < \tau < m < \tau' < M < \omega$ where $d(\tau) = d(\tau')$ is the absolute minimum value, $d(m)$ is a local maximum value, and $d(M)$ is the absolute maximum value. The circles $C(m) = \mu_{m}(C)$ and $C(M) = \mu_{M}(C)$ are orthogonal to $O$ and $C(\tau) =\mu_{\tau}(C)$ and $C(\tau')= \mu_{\tau'}(C)$ are centered in the first quadrant and are exchanged by the inversion $I_{O}$, which generates a periodic symmetry of the graph similar to the case articulated in the preceding paragraph.

\cref{fig:flowgraph_ellip} presents an example where $a= 5/3$ and the coaxial family $\mathcal{A}$ has foci $\pm y$, where $y = 0.5$. The three circles $C_{u}$, $C_{v}$ and $C_{w}$ are suppressed in the figure, and the respective values of $b$ and $c$ are $b = 5.846$ and $c=1.227$. In the close up of \cref{fig:flowgraph_ellip2}, the three local maximum circles (blue), the three absolute maximum circles (green), and the six enveloping circles are orthogonal to $O$ (not shown).


\subsection{Formul{\ae}}
The next order of business is to derive useful formul{\ae} for analyzing $d(t)$. Let $C(z,r)$ denote the circle in the complex plane centered at $z$ and of radius $r > 0$. The inversive distance between the circle $C(z,r)$ and its rotated cousin $\mathfrak{r}(C(z, r) ) = C(\mathfrak{r}(z), r)$ is
\begin{equation}\label{EQ:D(z,t)}
	h(z,r) = \frac{| z - \mathfrak{r}(z)|^{2} - 2 r^{2}}{2r^{2}} = \frac{1}{2}\frac{ \mid\sqrt{3} z - \boldsymbol{i}\mid^{2}}{r^{2}} -1.
\end{equation}
\begin{figure}
\includegraphics[width=0.35\textwidth]{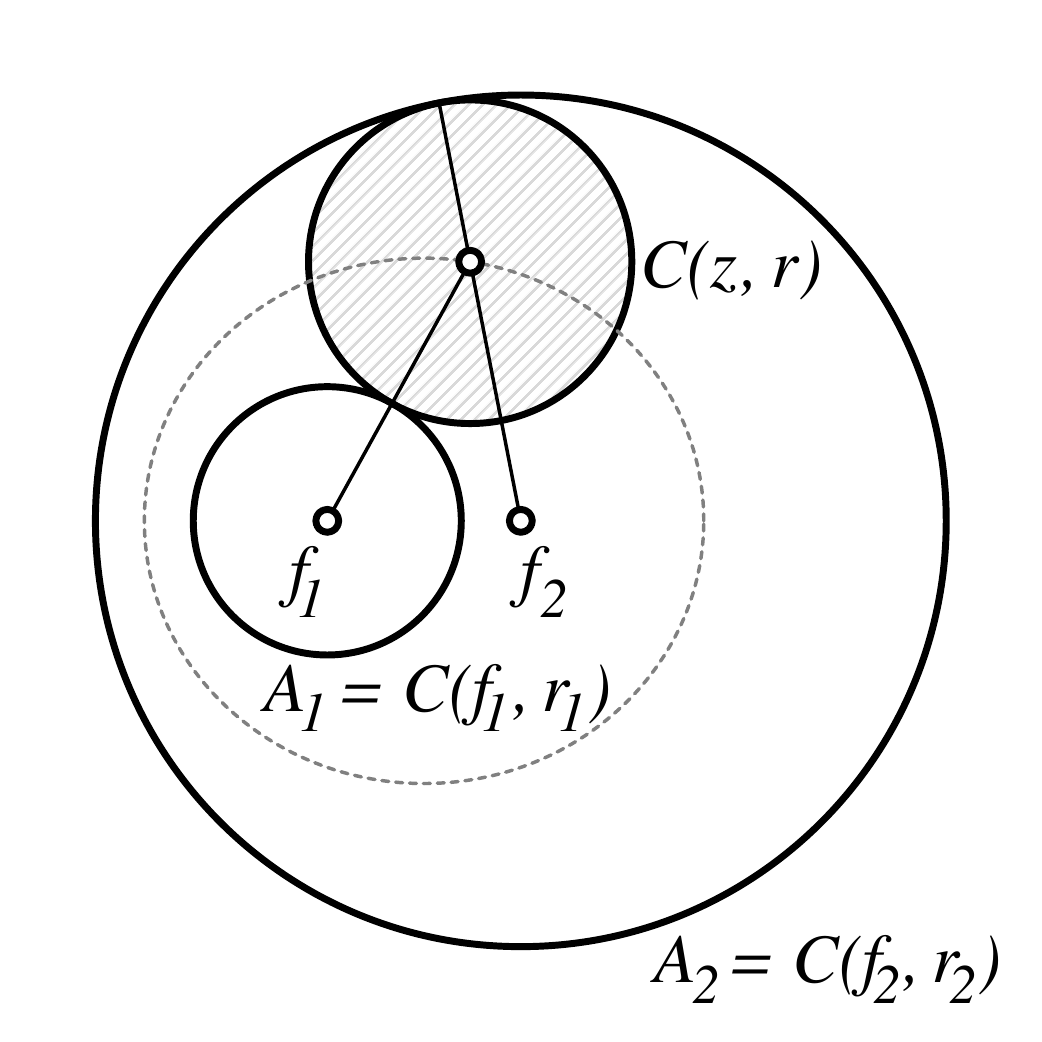}
\caption{The centers of $C(z, r)$ lie on the dotted ellipse with equation $|z-f_1|+|z-f_2|=r_1+r_2$.}
\label{fig:centers}
\end{figure}
\\
Let $f_{1}$ and $f_{2}$ be the respective centers of the two circles $A_{1}$ and $A_{2}$. A glance at \cref{fig:centers} should convince the reader that if $C(z,r)$ is externally tangent to $A_{1}$ and internally tangent to $A_{2}$, then $|z - f_{1} | = r_{1} + r$ and $|z - f_{2}| = r_{2} - r$, where $r_{1}$ and $r_{2}$ are the respective radii of the circles $A_{1}$ and $A_{2}$. It follows that
\begin{equation}\label{EQ:ellipse}
	|z - f_{1} | + |z - f_{2}|  = r_{1} + r_{2},
\end{equation}
implying that the centers $z = z(t)$ of the circles $\mu_{t}(C)$ lie on the ellipse described by Equation~\ref{EQ:ellipse}.  Using $|z - f_{1} | = r_{1} + r$ to eliminate $r$ from Equation~\ref{EQ:D(z,t)} gives the formula
\begin{equation}\label{EQ:h(z)}
	h(z) = \frac{1}{2}\frac{ \mid\sqrt{3} z - \boldsymbol{i}\mid^{2}}{\left(|z - f_{1} | -r_{1}\right)^{2}} -1
\end{equation}
for the inversive distance between a circle $C(z,r)$ externally tangent to $A_{1}$ and its rotated cousin $\mathfrak{r}(C(z,r))$. To find the extrema of the inversive distance function $d(t)$ of Equation~\ref{EQ:d(t)}, we need to find the extrema of $h(z)$ of Equation~\ref{EQ:h(z)} when $z$ is subject to the constraint of Equation~\ref{EQ:ellipse}. 

The ellipse of Equation~\ref{EQ:ellipse} has foci $f_{1}$ and $f_{2}$ with center $(f_{1} + f_{2})/2$, major radius $(r_{1} + r_{2})/2$, and minor radius $\frac{1}{2}\sqrt{(r_{1} + r_{2})^{2} - (f_{2} - f_{1})^{2}}$. The Cartesian equation for this ellipse is
\begin{equation}\label{EQ:EllipseEq}
	\frac{(2x - f_{1} - f_{2})^{2}}{(r_{1} + r_{2} )^{2}} + \frac{4y^{2}}{(r_{1} + r_{2})^{2} - (f_{2} - f_{1})^{2}} = 1.
\end{equation}

Equation~\ref{EQ:h(z)} has been used to compute the value of $d(\tau)$, the critical minimum value of $d(t)$ in the examples presented in this paper, and the values of $d(t) = d(t')$ for the Ma-Schlenker pairs. We now describe the derivations of the plots of \cref{fig:flowgraph_para,fig:flowgraph_ellip,fig:flowgraph_hyp}. In each case this is accomplished by mapping via a M\"obius transformation to a standard flow, applying the standard flow, and mapping back.

\begin{figure}
\includegraphics[width=0.75\textwidth]{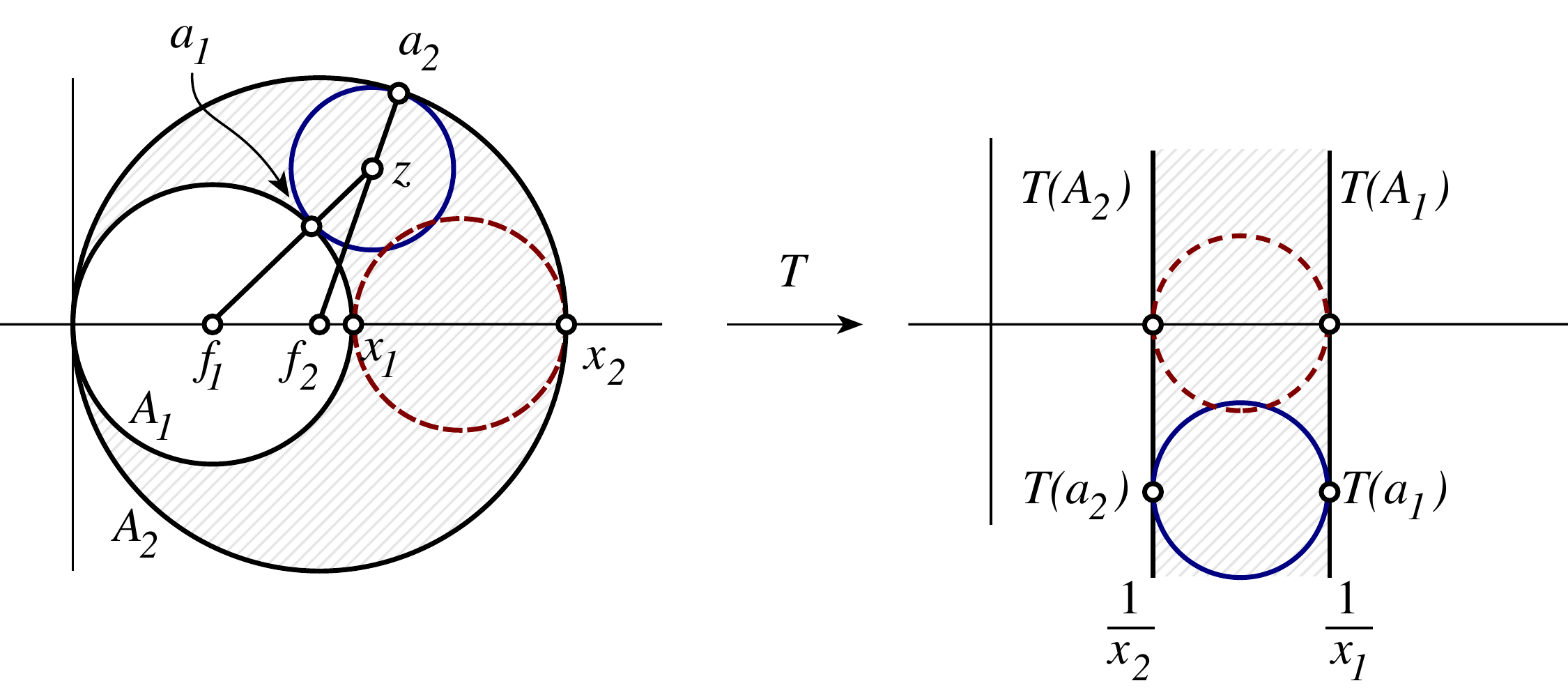}
\caption{The canonical M\"obius flow for the parabolic case. }
\label{fig:canpara}
\end{figure}
\subsubsection*{The parabolic case, \cref{fig:flowgraph_para}} Let $a_{1}= a_{1}(t)$ and $a_{2}=a_{2}(t)$ be the points of tangency of the circle $C(t) = \mu_{t}(C)$ with the respective enveloping circles $A_{1}$ and $A_{2}$, as in \cref{fig:canpara}. We will derive formul{\ae} for $a_{1}$ and $a_{2}$ in terms of parameters $1 < x_{1} < x_{2}$ and the variable $t$ and use these to write the center $z(t)$ and radius $r(t)$ of $C(t)$ for use in the function 
\begin{equation}\label{EQ:d=h(z,r)}
  d(t) = h (z(t),r(t))
\end{equation}
of Formula~\ref{EQ:D(z,t)}. The M\"obius transformation $T(z) = 1/z$, which is its own inverse, maps $A_{1}$ to the vertical line $x= 1/x_{1}$ and $A_{2}$ to the vertical line $x=1/x_{2}$. The flow  $\mu_{t}$ is obtained by conjugation of the standard unit speed flow $\nu_{t}(z) = z- t\boldsymbol{i}$ with $T$. The formul{\ae} for $a_{1}(t)$ and $a_{2}(t)$ become
\begin{equation*}
  a_{s}=a_{s}(t) = T\left(\frac{1- x_{s} t \boldsymbol{i}}{x_{s}}\right) = x_{s} \frac{1+ x_{s} t \boldsymbol{i}}{1 + x_{s}^{2}t^{2}}, \quad \text{for }s=1,2.
\end{equation*}
The radius $r(t)$  of $C(t)$ becomes
\begin{equation}\label{EQ:r(t)}
  r = r(t) = \frac{r_{1}r_{2}(a_{2}(t) - a_{1}(t))}{r_{2}(a_{1}(t)-f_{1})+ r_{1}(a_{2}(t) - f_{2})},
\end{equation}
and its center $z(t)$ becomes
\begin{equation}\label{EQ:z(t)}
  z= z(t) = f_{1} + (r_{1} + r(t) ) \frac{a_{1}(t) - f_{1}}{r_{1}}.
\end{equation}
\cref{fig:flowgraph_para} was generated using Formul{\ae}~\ref{EQ:r(t)} and \ref{EQ:z(t)} in Equation~\ref{EQ:d=h(z,r)} in \texttt{Mathematica}. We have posted \texttt{Mathematica} notebooks for generating visualizations of these constructions at \url{https://w3.cs.jmu.edu/bowersjc/page/circles/}.

\begin{figure}
\includegraphics[width=0.75\textwidth]{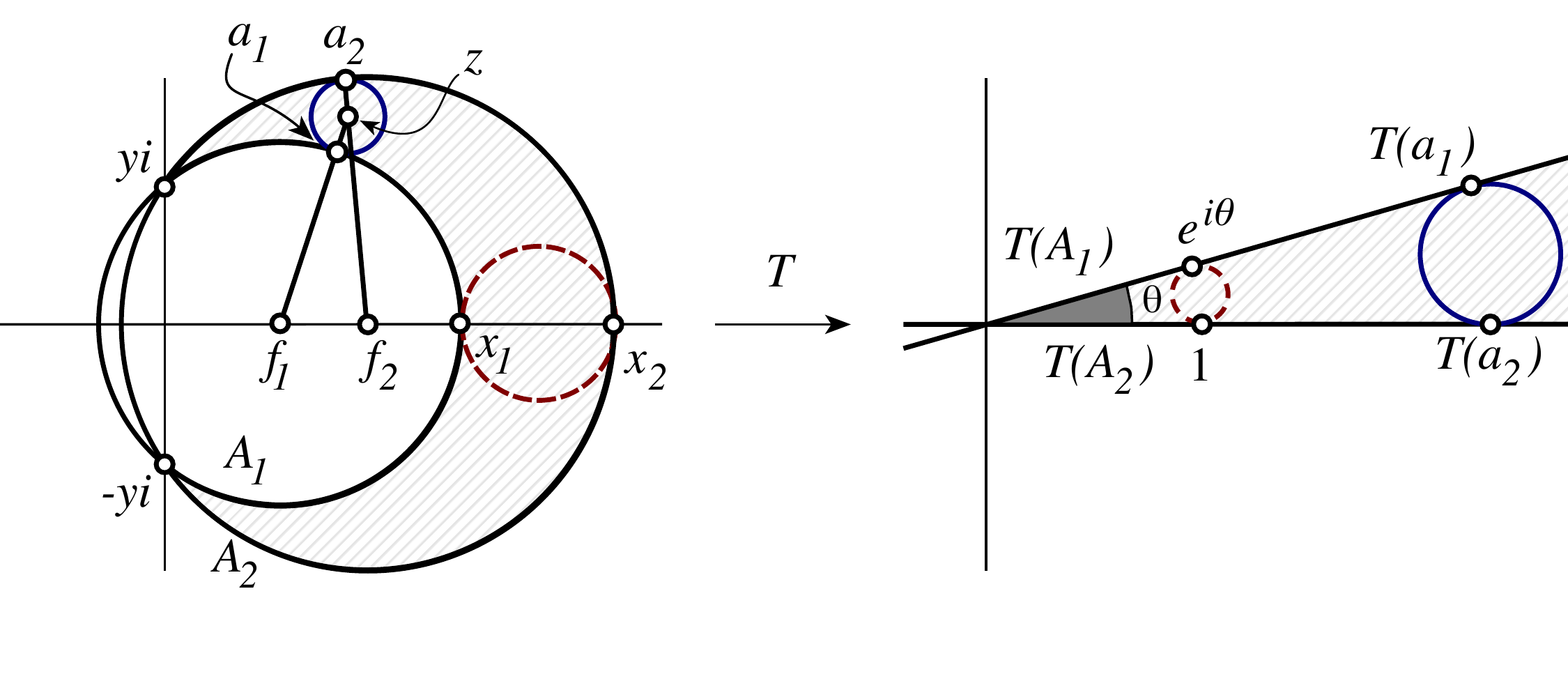}
\caption{The canonical M\"obius flow for the hyperbolic case. }
\label{fig:canhyp}
\end{figure}
\subsubsection*{The hyperbolic case, \cref{fig:flowgraph_hyp}} In this case the circles $C_{u}$ and $C_{v}$ meet at points $\pm y\boldsymbol{i}$ as in \cref{fig:canhyp}. The relationships among the parameters are that $y < 1/\sqrt{3}$, and for $s=1,2$, $f_{s}^{2} + y^{2} = r_{s}^{2}$ and $x_{s} = f_{s} + r_{s}$. The graphic is computed from Equations~\ref{EQ:d=h(z,r)}, \ref{EQ:r(t)}, and \ref{EQ:z(t)} with the only difference from the parabolic case that the formul{\ae} for $a_{1}$ and $a_{2}$ change. To obtain the correct formul{\ae}, we use the M\"obius transformation
\begin{equation}
  T(z) = \frac{x_{2} - y\boldsymbol{i}}{x_{2} + y\boldsymbol{i}}\; \frac{z+ y\boldsymbol{i}}{z- y \boldsymbol{i}}
\end{equation}
to map $-y\boldsymbol{i}$ to $0$, $y\boldsymbol{i}$ to $\infty$, and $x_{2}$ to $1$. The image of $A_{2}$ is the real axis and the image of $A_{1}$ is the line through the origin at angle 
\begin{equation}
  \theta = \cos^{-1} \left(\frac{f_{1}f_{2} + y^{2}}{r_{1}r_{2}}\right)
\end{equation}
up from the real axis. The standard unit speed flow is now $\nu_{t}(z) = e^{t} z$ and $a_{1}$ and $a_{2}$ are given as
\begin{equation}
  a_{1} = a_{1}(t) = T^{-1}\left(e^{t+ \theta\boldsymbol{i}}\right) \quad\text{and}\quad a_{2} = a_{2}(t) = T^{-1} (e^{t}), 
\end{equation}
where, setting $\kappa = (x_{2} - y\boldsymbol{i})/(x_{2} + y \boldsymbol{i})$, we may write $T$ and $T^{-1}$ as
\begin{equation}
  T(z) = \kappa\frac{z +  y \boldsymbol{i}}{z- y \boldsymbol{i}} \quad\text{and}\quad T^{-1} (z) = y\boldsymbol{i}\frac{z+ \kappa}{z-\kappa}.
\end{equation}

\subsubsection*{The elliptic case, \cref{fig:flowgraph_ellip}} In this case the circles $C_{u}$ and $C_{v}$ are part of a coaxial family with foci $\pm y$. The relationships among the various parameters are $0 < y < 1$, and for $s=1,2$, $f_{s} = (y^{2} + x_{s}^{2})/2x_{s}$ and $r_{s} = x_{s} - f_{s} = (y^{2}-x_{s}^{2})/2x_{s}$. Again, the graphic is computed from Equations~\ref{EQ:d=h(z,r)}, \ref{EQ:r(t)}, and \ref{EQ:z(t)} with $a_{1}$ and $a_{2}$ given by
\begin{equation}
  a_{1} = a_{1}(t) = T^{-1} \left(e^{t\boldsymbol{i} }\right) \quad \text{and} \quad a_{2} = a_{2}(t) = T^{-1} \left(e^{t\boldsymbol{i}} T(x_{2})\right).
\end{equation}
Here $T$ is the M\"obius transformation that takes $y$ to $0$, $-y$ to $\infty$, and $x_{1}$ to $1$. This takes the circle $A_{1}$ to the unit circle and $A_{2}$ to a circle centered at the origin of radius greater than $1$. The standard unit speed flow is the rotation flow $\nu_{t}(z) = e^{t\boldsymbol{i}} z$ and $T$ and $T^{-1}$ are given by
\begin{equation}
  T(z) = \kappa \frac{ z - y }{z + y} \quad\text{and}\quad T^{-1} (z) = -y\frac{z+ \kappa}{z-\kappa},
\end{equation}
where $\kappa = (x_{1} + y)/(x_{1} - y)$.

\quad


\begin{figure}
\includegraphics[width=0.3\textwidth]{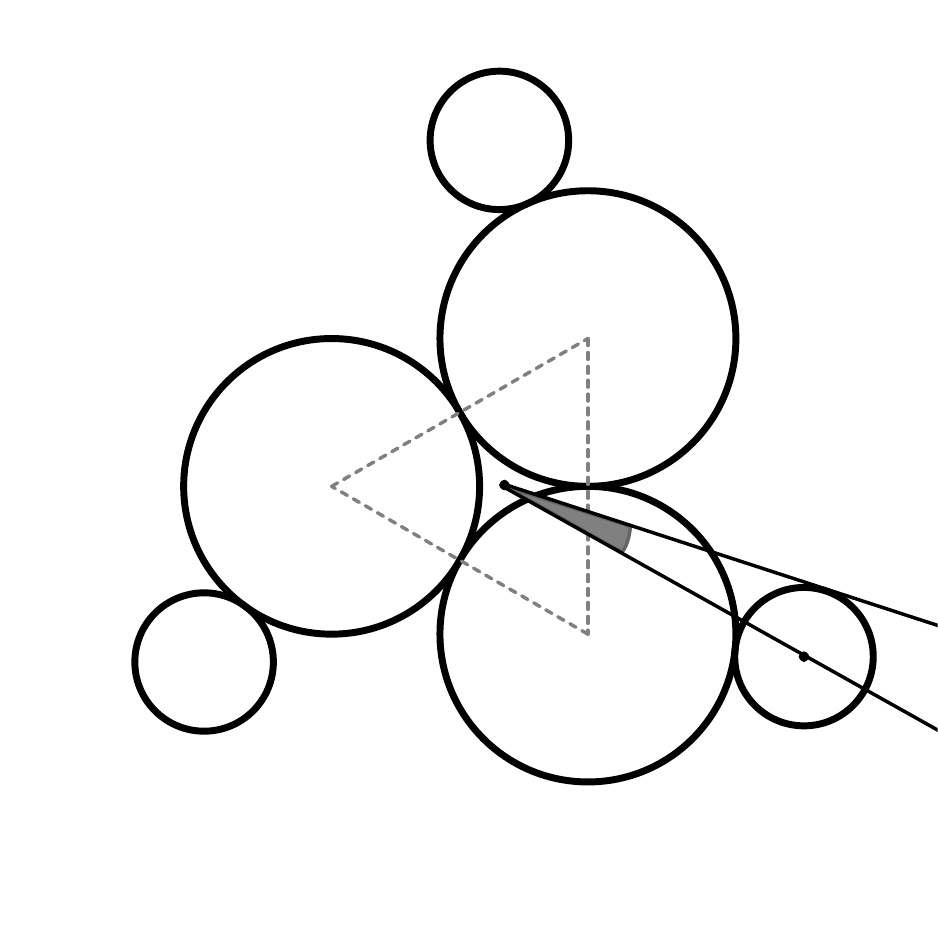}
\caption{The angle $\alpha$.}
\label{fig:figangle}
\end{figure}
\subsection{The angular equation}
One more geometric fact is useful. In the setting of the planar construction of Section~\ref{Section:Construction}, let $\alpha = \alpha(t)$ be the angle subtended by the ray from the incenter of the triangle $\Delta$ through the center of $C_{w'} = C_{w'}(t)$ and a ray from the incenter that is tangent to $C_{w'}$, as in \cref{fig:figangle}. It is easy to see that the angle $\alpha$ lies between $0$ and $\pi/2$, and an  algebraic manipulation determines that 
\begin{equation}\label{EQ:theta}
	d(t) = \langle C_{w'}, \mathfrak{r}(C_{w'}) \rangle = \frac{1}{2} + \frac{3}{2} \cot^{2} \alpha(t).
\end{equation}
Notice as $\alpha$ strictly increases from $0$ to $\pi/2$, the function $\cot^{2} \alpha$ strictly decreases from $\infty$ to $0$. Thus the comparison of two values of $d$ may be made simply by comparing the corresponding angles $\alpha$: 
\begin{equation}
	d(t) < d(t') \text{ if and only if } \alpha(t) > \alpha(t').
\end{equation}
From this observation one can make quick judgements of the correctness of many of our claims in the examples just by examining the graphics of the planar circle realizations of the figure, perhaps with a straight edge and protractor.

\section{The Construction on the Sphere}\label{Section:Examples}
\begin{figure}
\begin{subfigure}[t]{0.35\textwidth}
\includegraphics[width=\textwidth]{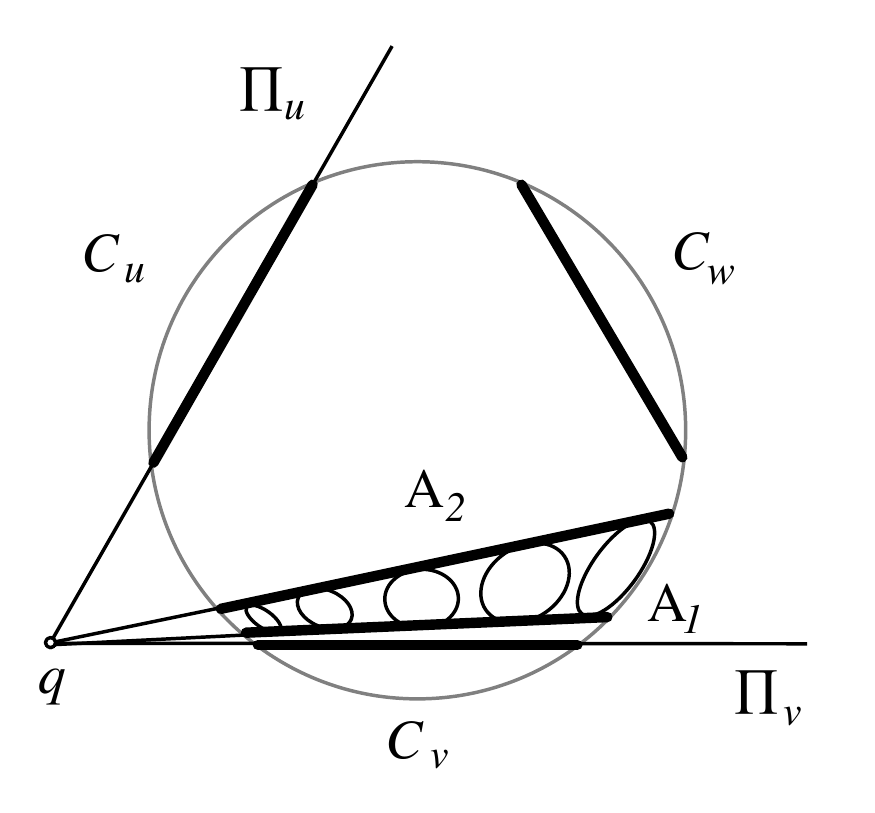}
\caption{Above the north pole.}
\label{fig:spheretop1}
\end{subfigure}
\quad
\begin{subfigure}[t]{0.30\textwidth}
\includegraphics[width=\textwidth]{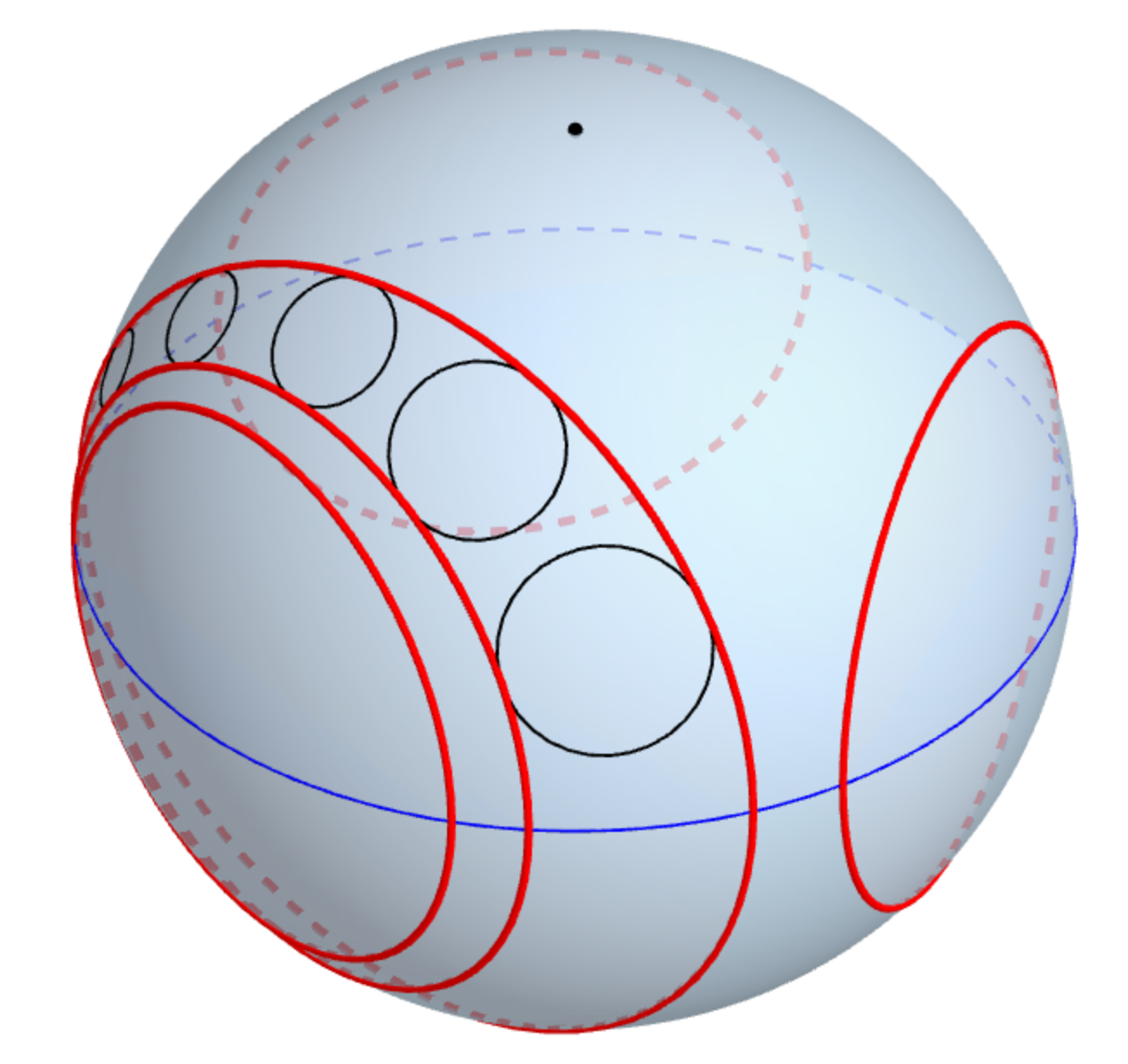}
\caption{A 3D view.}
\label{fig:sphereother}
\end{subfigure}
\caption{Two views of a state of the construction of a Ma-Schlenker realization with $C_{u}$, $C_{v}$, and $C_{w}$ centered on the equator.}
\label{fig:sphereconstruction}
\end{figure}

In this final section we present an alternative description of the construction of Ma-Schlenker realizations, this time directly on the $2$-sphere $\mathbb{S}^{2}$ realized as the unit sphere in $\mathbb{E}^{3} = \mathbb{C} \times \mathbb{R}$. In this normalization, the north pole is $\mathsf{n} = (0,1)$, the south pole is $\mathsf{s} = (0,-1)$, and the equatorial plane is identified with $\mathbb{C}$. The circles $C_{u}$, $C_{v}$, and $C_{w}$ in $\mathbb{S}^{2}$ are of equal radii whose centers are equally spaced on the equator. Projecting orthogonally along the north-south direction to the equatorial plane, the upper hemisphere projects to the unit disk and the three circles $C_{u}$, $C_{v}$, and $C_{w}$ project to symmetrically placed chords of the unit circle, as in \cref{fig:spheretop1}. The orthogonal projections of the typical circles $A_{1}$ and $A_{2}$ in the coaxial family $\mathcal{A}_{C_{u}, C_{v}}$ that form the envelope of the M\"obius-flowed circles $C_{w'}(t)$ are shown in \cref{fig:spheretop1}. The coaxial family $\mathcal{A}_{C_{u},C_{v}}$ is obtained as follows. Let $\Pi_{u}$ and $\Pi_{v}$ be the $2$-planes in $\mathbb{E}^{3}$ whose respective intersections with the $2$-sphere $\mathbb{S}^{2}$ are the circles $C_{u}$ and $C_{v}$. Let $\ell_{u,v}$ be the line $\Pi_{u} \cap \Pi_{v}$, a line that runs in the north-south direction. Then $\mathcal{A}_{C_{u},C_{v}}$ is precisely the collection of circles formed as intersections $\Pi \cap \mathbb{S}^{2}$ as $\Pi$ ranges over all $2$-planes in $\mathbb{E}^{3}$ that contain the line $\ell_{u,v}$. Notice that the circles in the coaxial family $\mathcal{A}_{C_{u},C_{v}}$ project orthogonally to the family of lines in $\mathbb{C}$ that pass through the common point $q$, the intersection of the vertical line $\ell_{u,v}$ with the equatorial plane $\mathbb{C}$. This explains the position of the projections of $A_{1}$ and $A_{2}$ in \cref{fig:spheretop1}. In this set up, local maxima for $d(t)$ occur when the three circles corresponding to $u'$, $v'$ and $w'$ are all centered on the equator. This occurs once when the flow $\mu_{t}$ is hyperbolic or parabolic, and periodically with two different local maximum values as the flow pushes the circle across the equator when elliptic. By shifting the origin, we may assume that a local maximum occurs at time $t=0$ so that the circles $C_{u'}(0)$, $C_{v'}(0)$, and $C_{w'}(0)$ are centered on the equator, and, without loss of generality, we may assume that the flow is oriented so that $C_{w'}(t)$ is centered in the upper hemisphere for initial positive values of $t$. Under these normalizations $d(t) = d(-t)$ and $C_{w'}(-t) = I_{O}(C_{w'}(t))$ for all $t$, where $O$ is the equator. The resulting realizations of the Ma-Schlenker octahedra of the form $\mathscr{O}(a,b,c,d)$ have order three rotational symmetry about the axis through the north and south poles. We should comment that a realization built in this manner will not be a circle packing. To obtain a packing these will need to be M\"obius flowed using a hyperbolic flow from the north toward the south pole.

\begin{figure}
\includegraphics[width=0.35\textwidth]{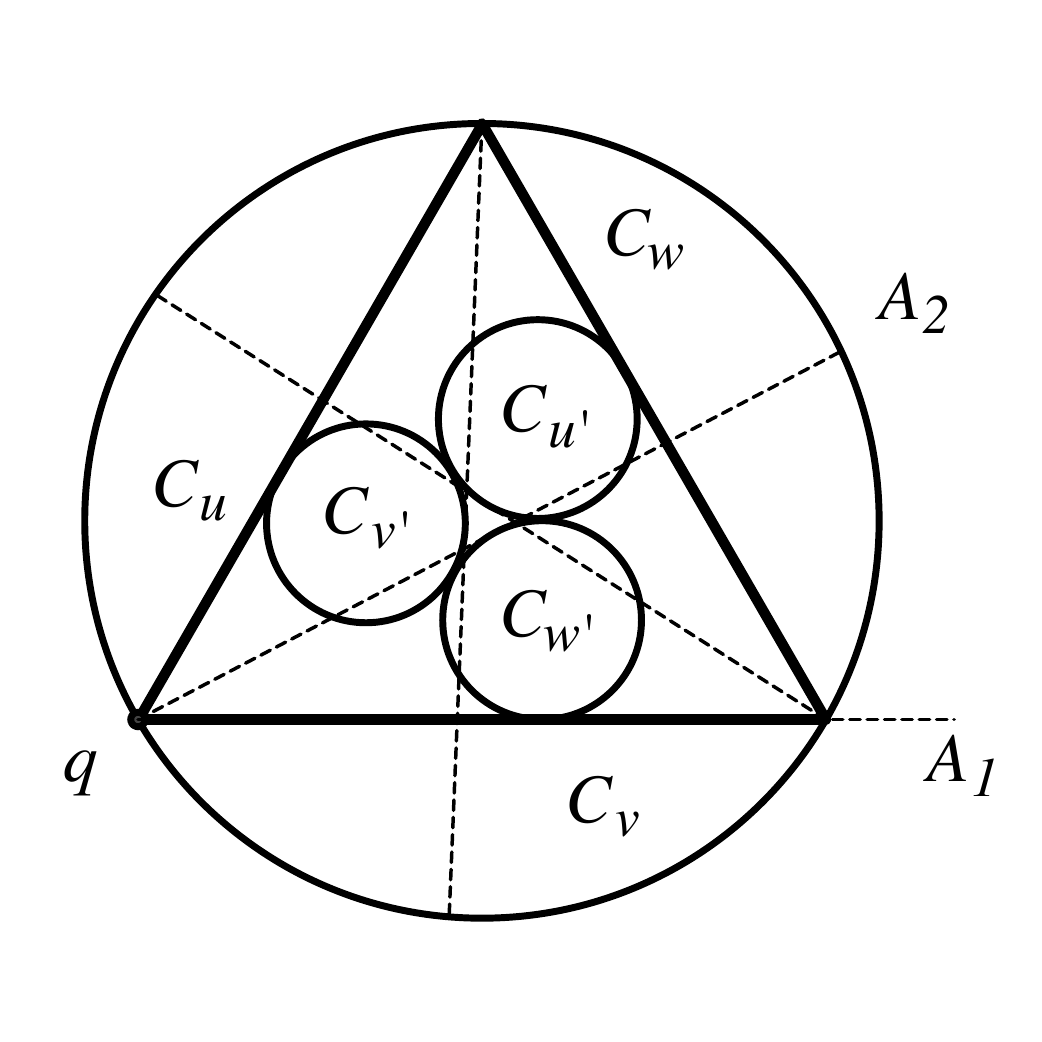}
\caption{The view of a critical circle realization for $\mathscr{O}(1,b,1,1)$ from above the north pole, with $C_{u}$, $C_{v}$, and $C_{w}$ centered on the equator.}
\label{fig:example41}
\end{figure}
As an example, we construct a critical circle packing for the Ma-Schlenker octahedron $\mathscr{O}(1, b,1,1)$ where $b>1$ and $d(\tau) = 1$. Choose $C_{u}$, $C_{v}$, and $C_{w}$ to have equal radii of $\pi/3$ so that the three circles are mutually tangent and project orthogonally to an equilateral triangle inscribed in the unit circle, as in \cref{fig:example41}. Choose $A_{1} = C_{v}$ and $A_{2}$ so that the circles $C_{u'}(t)$, $C_{v'}(t)$, and $C_{w'}(t)$ are disjoint except at two values $t = \pm\tau$ of the flow variable, when the circles are mutually tangent. Automatically, $d(t)$ must take on its isolated minimum values at $t=\pm \tau$ where $d(\pm \tau) =1$, with $d(t) >1$ for $t \neq \pm \tau$. To obtain a critical Ma-Schlenker \textit{c}-octahedron, apply the hyperbolic M\"obius flow from the north toward the south pole whose flow lines are the meridianal circles until the centers of $C_{u}$, $C_{v}$, and $C_{w}$ lie on a latitude $L$ and those of $C_{u'}(\tau)$, $C_{v'}(\tau)$, and $C_{w'}(\tau)$ lie on a latitude $L'$ of equal radii, $L$ centered on the south pole and $L'$ on the north. The resulting realization is a circle packing and therefore a critical Ma-Schlenker \textit{c}-octahedron.

\bibliographystyle{abbrv}
\bibliography{biblio}

\end{document}